\documentclass[12pt, reqno]{amsart}

\makeatletter
\def\theequation{\@arabic\c@equation}
\usepackage{amsmath,amsthm,amscd,amssymb,latexsym,upref}
\usepackage{versions}
\usepackage{enumerate}
\usepackage{graphicx}
\usepackage[usenames,dvipsnames] {color}

\def\d{\mathbb{D}}

\def\t{\mathbb{T}}

\def\be{\begin{equation}}
\def\ee{\end{equation}}

\def\s0{s_0}
\def\p0{p_0}
\let\phi\varphi
\let\epsilon\varepsilon

\newtheorem{theorem}{Theorem}[section]

\newtheorem{lemma}[theorem]{Lemma}

\numberwithin{equation}{section}

\newtheorem{lem}[theorem]{Lemma}
\newtheorem{prop}[theorem]{Proposition}
\newtheorem{thm}[theorem]{Theorem}

\newtheorem{definition}[theorem]{Definition}

\newtheorem{remark}[theorem]{Remark}

\newtheorem{fact*}[theorem]{Fact}
\newtheorem{defn}[theorem]{Definition}
\newtheorem{cor}[theorem]{Corollary}

\newtheorem{rem}[theorem]{Remark}
\newtheorem{exmp}[theorem]{Example}
\newtheorem{probl}[theorem]{Problem}

\newcommand\e{\mathrm{e}}
\newcommand\E{\mathcal{E}}

\newcommand\ii{\mathrm{i}}

\newcommand{\T}{\mathbb{T}}

\newcommand{\D}{\mathbb{D}}
\newcommand{\C}{\mathbb{C}}

\newcommand{\R}{\mathbb{R}}

\newcommand \z{\mathbb Z}

\newcommand{\ph}{\varphi}
\newcommand\al{\alpha}

\newcommand\la{\lambda}

\newcommand\si{\sigma}
\newcommand\ta{\theta}

\newcommand\beq{\begin{equation}}

\newcommand\eeq{\end{equation}}
\newcommand\ov{\overline}

\newcommand\bbm{\begin{bmatrix}}
\newcommand\ebm{\end{bmatrix}}
\newcommand\bpm{\begin{pmatrix}}
\newcommand\epm{\end{pmatrix}}

\let\phi\varphi
\numberwithin{equation}{section}

\begin{document}
\title[Interpolation by tetra-inner maps]{Interpolation by holomorphic maps from the disc to the tetrablock}

\author{Hadi O. Alshammari and Zinaida A. Lykova}
\date{6th January 2021}

\begin{abstract} 
The {\em tetrablock} is the set
	$$
	\mathcal{E}=\{x \in \mathbb{C}^3: \quad 1-x_1z-x_2w+x_3z w  \neq   0 \quad  whenever \quad |z|\leq 1, |w|\leq 1\}.
	$$
The closure of $\mathcal{E}$ is denoted by $\overline{\mathcal{E}}$. A {\em tetra-inner function} is an analytic map $x$ from the unit disc $ \mathbb{D} $  to $\overline{\mathcal{E}}$ such that, for almost all points $\la$ of the unit circle $ \mathbb{T}$,
\[
\lim_{r\uparrow 1} x(r\la) \mbox{ exists and lies in } b \overline{\mathcal{E}},
\]
where $b \overline{\mathcal{E}}$ denotes the distinguished boundary  of $\overline{\mathcal{E}}$. There is a natural notion of {\em degree} of a  rational tetra-inner function $ x$; it is simply the topological degree of the  continuous map $ x|_\mathbb{T} $ from  $ \mathbb{T} $ to $ b \overline{\mathcal{E}} $.
	
	In this paper we give a prescription for the construction of a general rational tetra-inner function of degree $n$.  The prescription exploits a known construction of the finite  Blaschke products of given degree which satisfy some interpolation conditions with the aid of  a Pick matrix formed from the interpolation data. It is known  that if $x= (x_1, x_2, x_3)$ is a rational tetra-inner function of degree $n$, then  $x_1 x_2 - x_3$ either is identically  $0$ or  has precisely $n$ zeros in the closed unit disc $\overline{\mathbb{D}}$, counted with multiplicity.  
	It turns out that a natural choice of data for the construction of a rational tetra-inner function $x= (x_1, x_2, x_3)$	consists of the points in $\overline{\mathbb{D}}$ for which $x_1 x_2 - x_3=0$ and the values of $x$ at these points.
\end{abstract}



\thanks{The first author  was supported by the Government of Saudi Arabia. The second author was partially supported by the UK Engineering and Physical Sciences Research Council grant  EP/N03242X/1. }

\maketitle
\tableofcontents
\section{Introduction}\label{Nintro}

In this paper we present an algorithm for the construction of a general rational inner function from  $\D$ to the tetrablock. The algorithm is based on a known solution of the Nevanlinna-Pick interpolation problem on $\D$. Different versions of the Nevanlinna-Pick interpolation problem have been studied by many authors, beginning with G. Pick  in 1916 \cite{Pick1916} and continuing
with R. Nevanlinna in 1922 \cite{Nev1922}, and they still attract interest, because they are natural questions in function theory and because of their applications to engineering, particular electric networks and control theory, see \cite{Ball1983} for some references.
We should mention particularly papers of J. A. Ball and J. W. Helton \cite{BH},
D. Sarason \cite{Sar}, D. R. Georgijevi\'c \cite{Geo98}, G.-N. Chen and Y.-J. Hu \cite{ChH} and V. Bolotnikov and A. Kheifets \cite{BolKh08}, and the books of J. A. Ball, I. C. Gohberg and L. Rodman \cite{bgr}, and of V. Bolotnikov and H. Dym \cite{BD}.  There are many further papers on this interesting topic and applications (see, for example, \cite{CS,SeWe}).

The closed tetrablock $\overline{\mathcal{E}}$ is the set in $\mathbb{C}^{3}$ defined by
$$
\overline{\mathcal{E}} = \{(x_1,x_2, x_3) \in \mathbb{C}^3: \ 1-x_1z-x_2w+x_3zw \neq  0 \ \text{whenever} \ |z| < 1, |w| < 1\}.
$$
The original motivation for the study of $\overline{\mathcal{E}}$ was an attempt to solve a  $\mu$-synthesis problem \cite{AWY}, which is itself motivated by basic unsolved problems in $H^\infty$ control theory \cite{D,dullerud}. The tetrablock has attracted considerable interest in recent years. It has interesting complex geometry \cite{AWY,EKZ,Y1,KZ2015,KZ2016}, rich function theory \cite{BLY,AlsLyk,OZY} and associated operator theory \cite{TB}.
 The solvability of the $\mu$-synthesis problem connected to $\E$  can be expressed in terms of the existence of rational inner functions from the open unit disc $\D$ in the complex plane $\C$ to the closure of $\E$ \cite{BLY}. 

Observe that if $x:\D\to \ov{\mathcal{E}}$ is analytic then, by Fatou's Theorem, for almost all $\la\in\T$ with respect to Lebesgue measure, the radial limit
\[
\lim_{r\uparrow 1} x(r\la) 
\]
exists.   We say that an analytic map $x:\D\to \ov{\mathcal{E}}$ is a {\em tetra-inner function}, or alternatively, an {\em $\overline{\mathcal{E}}$-inner function} if, for almost all $\la \in\T$ with respect to Lebesgue measure,
\[
\lim_{r\uparrow 1} x(r\la) 
\]
lies in the distinguished boundary $b \overline{\mathcal{E}}$ of $\overline{\mathcal{E}}$.   The distinguished boundary $b \overline{\mathcal{E}}$ of 
$\overline{\mathcal{E}}$  is homeomorphic to the solid torus $\overline{\mathbb{D}} \times \T$,  which has a boundary \cite{AWY}.  The
$ \overline{\mathcal{E}}$-inner functions constitute a natural analogue (in the context of the tetrablock) of the {\em inner functions} introduced by A. Beurling \cite{ABB}, which play an important role in the function theory of the unit disc and in operator theory \cite{CGP}.  

A basic question about rational inner functions $\varphi$ from $\mathbb{D}$ to $\overline{\mathbb{D}}$ was  studied by W. Blaschke \cite{blaschke}.  Specifically, he obtained (inter alia) a formula for the general  rational inner function $\varphi$ of degree $n$ in terms of its zeros.  Indeed, by the Argument Principle, any rational inner function $\varphi$ of degree $n$ has exactly $n$ zeros in $\mathbb{D}$, counted with multiplicity. From this fact one can see that  $\varphi$  is a ``finite Blaschke product", having the form
\be\label{Bla}
\varphi(\lambda) = c  {\displaystyle \prod_{j=1}^{n}} \frac { \lambda-\alpha_j} {1-\overline{\alpha_j}\lambda}\ 
\ee
for some unimodular constant $c$ and some $\alpha_1, ... , \alpha_n\in\D$.  The $\al_j$ are the zeros of $\varphi$.  It is evident from equation \eqref{Bla} that $\varphi$ extends to a continuously differentiable function on $\overline{\D}$, given by the same formula.
In this paper our aim is to write down a formula analogous to equation \eqref{Bla} for the general rational $\overline{\mathcal{E}}$-inner function of degree $n$.   The first question that arises is: what data should replace the $\al_j$, the zeros of $\ph$?  We have found that an effective choice is the set of {\em royal nodes} of the tetra-inner function, which we shall now define.
It was shown in \cite{AlsLyk} that if  $x=(x_1, x_2, x_3)$ is a rational $\overline{\mathcal{E}}$-inner function of degree $n$ then $x_3 -x_1 x_2$
either is identically $0$ or  has exactly $n$ zeros in the closed unit disc $\overline{\d}$, counted with an appropriate notion of multiplicity.
Here, the {\em degree} of a rational $\overline{\mathcal{E}}$-inner function $x$ is defined to be the topological degree of the restriction of $x$ that maps $\t$ continuously to $b \overline{\mathcal{E}}$.  Since $b \overline{\mathcal{E}}$ is homeomorphic to the solid torus $\overline{\mathbb{D}} \times \T$,  which is homotopic to $\T$, the fundamental group $\pi_1(b \overline{\mathcal{E}})$ is $\z$, and so the degree of $x$ is an integer; it will be denoted by $\deg(x)$.

The variety 
$$
 \mathcal{R}_{\overline{\mathcal{E}}} =   \{ (x_1, x_2, x_3) \in \overline{\mathcal{E}} : x_3 = x_1 x_2        \}
 $$
has an important role in the function theory of ${\mathcal{{E}}}$; it is called the {\em royal variety of $\overline{\mathcal{E}}$}. For any rational tetra-inner function $x = (x_1, x_2, x_3)$, the zeros of $x_1 x_2 - x_3$ in $\overline{\mathbb{D}}$ are the points $\lambda \in \mathbb{\overline{D}}$ such that $x(\lambda) \in \mathcal{R}_{\bar{\mathcal{E}}} $. We call these points  the {\em royal nodes of $x$}. If $\sigma \in \overline{\mathbb{D}}$ is a royal node of $x$, so that $x(\sigma) = (\eta, \tilde{\eta}, \eta \tilde{\eta})$ for some $\eta, \tilde{\eta} \in \overline{\mathbb{D}}$, then we call $\eta, \tilde{\eta}$ the {\em royal values of $x$} corresponding to the royal node $\sigma$ of $x$. In this paper, in Theorem \ref{main-theorem}, we give a prescription for the construction of a general 
 rational tetra-inner function of degree $n$ in terms of its royal nodes and royal values.  We shall make use of  a known solution of an interpolation problem for finite Blaschke products. 

The $\overline{\mathcal{E}}$-inner functions $x$ such that $x(\d)\subseteq \mathcal{R}_{\overline{\mathcal{E}}}$ are simply the functions of the form  $(\ph_1,\ph_2, \ph_1\ph_2)$ where $\ph_1, \ph_2$ are  inner functions.  These $\overline{\mathcal{E}}$-inner functions behave differently from the others, and we shall often specifically exclude them from consideration.

To describe our main results  on the construction of a general rational tetra-inner function we need to recall some definitions and results on the Blaschke interpolation problem. 
\begin{defn} {\em \cite[Definition 1.2] {ALY3}} \label{BlaschkeInterpolationData}
	Let $n \geq 1 $ and $ 0 \leq k \leq n$. By {\em Blaschke interpolation data of type $(n,k)$}  we mean a triple $(\sigma , \eta , \rho ) $ where 
	\begin{enumerate}
		\item[\em (i)] $ \sigma = (\sigma_1, \sigma_2, ... , \sigma_n)$ is an $n$-tuple of distinct points of \ $\overline{\mathbb{D}}$ such that $ \sigma_j \in \mathbb{T} $ for $ j = 1 , ... , k $ and $  \sigma_j \in \mathbb{D} $ for $ j = k + 1, ... , n;$
		\item[\em (ii)] $ \eta = (\eta_1, \eta_2, ..., \eta_n)$ where $ \eta_j \in \mathbb{T}$ for $ j = 1, ... ,k $ and $ \eta_j \in \mathbb{D} $ for $ j = k + 1, ... , n;$
		\item[\em (iii)] $ \rho = (\rho_1, \rho_2, ..., \rho_k)$ where $ \rho_j \  > \ 0 $ for $ j = 1, ... ,k $.
	\end{enumerate}
\end{defn}

\begin{definition}\label{phasar-derivative} {\rm \cite[Definition 7.5]{ALY2}~}
 For any differentiable function $f: \T \to \C \setminus \{0\}$
 the {\em phasar derivative} of $f$ at $z= \e^{\ii \ta} \in \T$
is  the derivative with respect to $\ta$ of the argument of $f(\e^{\ii \ta})$ at $z$; we denote it by $Af(z)$.
\end{definition}
Thus, if $f(\e^{\ii \ta}) = R(\ta) \e^{\ii g(\ta)}$ is differentiable,   where $g(\ta) \in\R$ and $R(\ta) >0$, then $g$ is differentiable on $[0,2\pi)$ and the phasar derivative of $f$ at $z= \e^{\ii \ta} \in \T$ is equal to
\begin{equation}\label{defA}
Af(\e^{\ii \ta}) = \frac{d}{d \ta} \arg f(\e^{\ii \ta}) = g'(\ta).
\end{equation}
We shall find it useful in the sequel. 

\begin{probl} \label{The Blaschke interpolation problem} \textbf{\em{(The Blaschke interpolation problem)}} 
	For given Blaschke interpolation data $(\sigma , \eta , \rho )$ of type $(n,k)$, find all  rational inner functions  $ \varphi $ on $ \mathbb{D}$  of degree $n$ with the properties 
	\begin{equation}
	\varphi(\sigma_j)=\eta_j\quad \text{ for} \ j=1,...,n
	\end{equation}
	and
	\begin{equation}\label{eq6.4}
	A\varphi(\sigma_j)=\rho_j \quad \text{for} \ j=1,...,k.
	\end{equation}
\end{probl}

Problem \ref{The Blaschke interpolation problem} has been analysed by several authors \cite{RR,Sar,BH}.  In the absence of the tangential conditions
\eqref{eq6.4} the problem would be ill-posed, in that the solvability of the problem would depend only on the interpolation conditions at nodes in $\d$, and the conditions at $\si_1,\ldots,\si_k$ would be irrelevant.  With the conditions \eqref{eq6.4}, however, the problem has a pleasing solution.
The existence of a solution of the Blaschke interpolation problem can be characterized in terms of an associated ``Pick matrix", and all solutions $\varphi$ are parametrized by a linear fractional expression in terms of a parameter $\zeta \in \mathbb{T}$. There are polynomials $a, b, c$ and $d$ of degree at most $n$ such that the general solution of Problem \ref{The Blaschke interpolation problem} is
$$
\varphi = \frac{a \zeta +b}{c \zeta +d}
$$ 
where the parameter $\zeta$ ranges over a cofinite subset of $\mathbb{T}$.
 The polynomials $a,b,c$ and $d$ are unique subject to a certain normalization.

\begin{defn}\label{prop3.10}{\em \cite[Definition 3.10] {ALY3}}
	Let $ (\sigma, \eta, \rho) $ be Blaschke interpolation data of type $(n,k)$. Suppose that {\em Problem \ref{The Blaschke interpolation problem} } is solvable. We say that
	$$
	\varphi = \frac { a\zeta + b} { c\zeta + d}
	$$
	is a {\em normalized linear fractional parametrization of the solutions }of {\em Problem \ref{The Blaschke interpolation problem} } if
	\begin{enumerate}
		\item[\em (i)] $a, b, c, d$ are polynomials of degree at most $n$;
		\item[\em (ii)] for all but at most $k$ values of $\zeta \in \mathbb{T} $, the function
		\begin{equation} \label{equ in def 3.10 new}
		\varphi(\lambda) = \frac { a(\lambda)\zeta + b(\lambda)} { c(\lambda)\zeta + d(\lambda) }
		\end{equation}
		is a solution of {\em Problem \ref{The Blaschke interpolation problem}};
		\item[\em (iii)] for some point $ \tau \in \mathbb{T} \ \backslash  \{ \sigma_1, ... , \sigma_k \} $, 
		$$
		\begin{bmatrix}
		a(\tau)       & b(\tau) \\
		c(\tau)       & d(\tau) \\
		\end{bmatrix}
		=
		\begin{bmatrix}
		1 & 0  \\
		0 & 1 \\
		\end{bmatrix};
		$$
		\item[\em (iv)] every solution $ \varphi $ of {\em Problem \ref{The Blaschke interpolation problem} } has the form {\em (\ref{equ in def 3.10 new})} for some $\zeta \in \mathbb{T}$. 		
	\end{enumerate}
\end{defn}

The data for the construction of a rational tetra-inner function $x$ consists of the royal nodes and royal values of $x$.

\begin{defn} \label{royaltetradata}
	Let $ n \geqslant 1 $, and $ 0\leq k \leq n $. By {\em royal tetra-interpolation data} of type $(n,k)$  we mean a four-tuple  $(\sigma, \eta, \tilde {\eta}, \rho) $ where 
	\begin{enumerate}
		\item[\em (i)] $ \sigma = (\sigma_1, \sigma_2, ... , \sigma_n)$ is an n-tuple of distinct points such that $ \sigma_j \in \mathbb{T} $ for $ j = 1 , ... , k $ and $  \sigma_j \in \mathbb{D} $ for $ j = k + 1, ... , n;$
		\item[\em (ii)] $ \eta = (\eta_1, \eta_2, ..., \eta_n)$ where $ \eta_j \in \mathbb{T}$ for $ j = 1, ... ,k $ and $ \eta_j \in \mathbb{D} $ for $ j = k + 1, ... , n;$
		\item[\em (iii)] $ \tilde {\eta} = ( \tilde {\eta_1}, \tilde {\eta_2}, ... , \tilde {\eta_n} ) $  where $ \tilde {\eta_j} \in \mathbb{T}$ \ for $ j = 1, ... , k $ and $ \tilde {\eta_j} \in \mathbb{D} $ \ for $ j = k+1, ... , n $. 
			\item[\em (iv)] $ \rho = (\rho_1, \rho_2, ..., \rho_k)$ where $ \rho_j \  > \ 0 $ for $ j = 1, ... ,k $.
	\end{enumerate}
\end{defn}

\begin{probl} \label{The royal tetra-interpolation problem} \textbf{\em (The royal tetra-interpolation problem)} Given royal tetra-interpolation 
data $(\sigma , \eta , \tilde {\eta}, \rho ) $ of type $(n,k)$,  find all rational $\overline{\mathcal{E}}$-inner functions $x=(x_1,x_2, x_3)$ of degree $n$ such that 
	$$
	x(\sigma_j)=(\eta_j, {\tilde {\eta}}_j, \eta_j {\tilde {\eta}}_j  ) \quad \text{for} \ j=1,...,n,
	$$
	and
	$$
	Ax_1(\sigma_j)=\rho_j \quad \text{for} \  j=1,...,k.
	$$	
\end{probl}

The connection between Problems \ref{The royal tetra-interpolation problem} and \ref{The Blaschke interpolation problem} can be described with the aid of a certain $1$-parameter family of rational functions  on $ \overline{\mathcal{E}}$ that is parametrized by the unit circle $\T$. These functions play a central role in the function theory of $ \overline{\mathcal{E}}$ (see \cite{AWY,BLY}).  They are defined, for $\omega\in\T$, by 
\be\label{eq1.26_I}
\Psi_{\omega}(x)  \ = \ \frac {  x_3 \omega - x_1 }  { x_2  \omega - 1  }  \ \  \textnormal{when} \quad x_{2} \omega - 1 \neq 0.
\ee
$\Psi_\omega$ is holomorphic on $ \overline{\mathcal{E}}$, except at points $ x \in  \overline{\mathcal{E}}$ where
$x_{2} \omega - 1 = 0$, and maps any point of $ \overline{\mathcal{E}}$  at which it is defined  into $\overline{\d}$.

The main theorem of this paper is the following.  

\begin{thm} \label{main-theorem}
	For royal tetra-interpolation data $(\sigma, \eta, \tilde{\eta}, \rho)$ of type $(n,k)$ the following two statements are equivalent:
	\begin{enumerate}
		\item[\em (i)] The royal tetra-interpolation problem {\em (Problem \ref{The royal tetra-interpolation problem})} with data $(\sigma, \eta, \tilde{\eta}, \rho)$  admits a rational $\overline{\mathcal{E}}$-inner interpolating function $x$ such that $x(\overline{\mathbb{D}}) \nsubseteq \mathcal{R}_{\overline{\mathcal{E}}}; $ 
		\item[\em (ii)] The Blaschke interpolation problem {\em (Problem \ref{The Blaschke interpolation problem})} with data $(\sigma, \eta, \rho)$ of type $(n,k)$ is solvable and there exist $x_1^{\circ}, x_2^{\circ}, x_3^{\circ} \in \mathbb{C}$ such that 
		$$
		| x_3^{\circ} | = 1, \ \ \ | x_1^{\circ} | < 1, \ \ \  | x_2^{\circ} | < 1, \ \	x_1^{\circ} = \overline{x_2^{\circ}} x_3^{\circ},
		$$ 
		and 
		$$
		\frac{x_3^{\circ} c(\sigma_j) + x_2^{\circ} d(\sigma_j)}{x_1^{\circ} c(\sigma_j) + d(\sigma_j)} = \tilde{\eta_j}  \quad \mbox{for} \ j = 1, ..., n,
		$$
		where $a, b, c$ and $d$ are the polynomials in the normalized parametrization $\varphi = \displaystyle{\frac{a \zeta + b}{c \zeta +d}}$ of the solution of Problem {\em \ref{The Blaschke interpolation problem}}.
	\end{enumerate}
\end{thm}
 The theorem follows from Theorems \ref{4.4Tetra} and \ref{ReverseThm}. 
The proofs of these theorems are given in Section \ref{FromRoyaltoBlash} and Section \ref{FromBlashtoRoyal} respectively. 
Theorem \ref{ReverseThm} provides us with a formula for a solution $x$ of Problem \ref{The royal tetra-interpolation problem} in terms of $s_0,p_0,a,b,c$ and $d$.  The formula is in terms of the polynomials $a,b,c$ and $d$ computed in \cite[Theorem 3.9]{ALY3}
 (see Remark \ref{explicit-solution}).  In this way we derive an explicit solution of  Problem \ref{The royal tetra-interpolation problem}. 

 {\bf Theorem \ref{ReverseThm}}. 
	{\it Let $(\sigma, \eta,  \rho)$ be Blaschke interpolation data of type $(n,k)$, and let $(\sigma, \eta, \tilde{\eta}, \rho)$ be royal tetra-interpolation data, where $\tilde{\eta}_j \in \mathbb{T}$, for $j = 1, ..., k$, and $\tilde{\eta}_j \in \mathbb{D}$, for $j = k+1, ..., n$. Suppose that the Blaschke interpolation problem (Problem {\em \ref{The Blaschke interpolation problem}}) with these data is solvable and the solutions $\varphi$ of Problem {\em \ref{The Blaschke interpolation problem}} have normalized parametrization 
	$$
	\varphi = \frac { a\zeta + b} { c\zeta + d}.
	$$
	Suppose that there exist scalars $x_1^{\circ}, x_2^{\circ}, x_3^{\circ} \in \mathbb{C}$ such that 
	$$
	| x_3^{\circ} | = 1, \ \ \ | x_1^{\circ} | < 1, \ \ \  | x_2^{\circ} | < 1, \ \	x_1^{\circ} = \overline{x_2^{\circ}} x_3^{\circ},
	$$ 
	and 
	\begin{equation} 
	\frac{x_3^{\circ} c(\sigma_j) + x_2^{\circ} d(\sigma_j)}{x_1^{\circ} c(\sigma_j) + d(\sigma_j)} = \tilde{\eta_j}  \quad \text{for} \  j = 1, ..., n. 
	\end{equation}
	Then there exists a rational tetra-inner function $x = (x_1, x_2, x_3)$  given by, 
	\begin{equation} 
	x_1(\lambda) = \frac { x_1^{\circ} a(\lambda) + b(\lambda) } { x_1^{\circ} c(\lambda) + d(\lambda) }
	\end{equation}
	\begin{equation} 
	x_2(\lambda) = \frac{x_3^{\circ} c(\lambda) + x_2^{\circ} d(\lambda) } { x_1^{\circ} c(\lambda) + d(\lambda) }
	\end{equation}
	\begin{equation} 
	x_3(\lambda) = \frac { x_2^{\circ} b(\lambda) + x_3^{\circ} a(\lambda)} { x_1^{\circ} c(\lambda) + d(\lambda) },
	\end{equation} for $\lambda \in \mathbb{D}$, such that
	\begin{enumerate}
		\item[\em (i)] $x$ is a solution of the royal tetra-interpolation problem with the data $(\sigma, \eta, \tilde{\eta}, \rho)$, that is, 
		$$x(\sigma_j) = (\eta_j, \tilde{\eta_j}, \eta_j \tilde{\eta_j}) \ \text{for} \ j = 1, ... ,n,$$
		and 		
		$$Ax_1(\sigma_j) = \rho_j \ \  \text{for} \  j = 1, ... , k,$$
		\item[\em (ii)] for all but finitely many $\omega \in \mathbb{T}$, the function $\Psi_{\omega} \circ x $ is a solution of Problem {\em \ref{The Blaschke interpolation problem}}. 
	\end{enumerate}}				

The solution sets of the royal tetra-interpolation problem and the corresponding Blaschke interpolation problem admit an explicit connection in terms of the  functions $\Psi_{\omega}$. 

{\bf Corollary \ref{LastCor}}. 
{\it	Let $(\sigma, \eta, \rho)$ be Blaschke interpolation data of type $(n,k)$. Suppose that $x$ is a solution of the {\em Problem \ref{The royal tetra-interpolation problem}} with data $(\sigma, \eta, \tilde{\eta}, \rho)$ for some $\tilde{\eta_{j}} \in \overline{\mathbb{D}}, j = 1, ..., n,$ and that $x(\mathbb{D}) \not\subseteq \mathcal{R_{\overline{\mathcal{E}}}}$. For all $\omega \in \mathbb{T} \backslash \{ \overline{\tilde{\eta_1}}, ... , \overline{{\tilde {\eta}}_k}\}$, the function $\varphi = \Psi_{\omega} \circ x$ is a solution of {\em Problem \ref{The Blaschke interpolation problem}} with Blaschke interpolation data $(\sigma, \eta, \rho)$. Conversely, if $\varphi$ is a solution  of the Blaschke interpolation problem with data $(\sigma, \eta, \rho)$, then there exists $\omega \in \mathbb{T}$ such that $\varphi = \Psi_{\omega} \circ x$ .		
}

In \cite{AlsLyk} there is a construction of the general rational $\overline{\mathcal{E}}$-inner function $x= (x_1,x_2, x_3)$ of degree $n$, in terms of different data, namely, the royal nodes of $x$ and the zeros of $x_1$ and $x_2$.  A major step in the construction in \cite{AlsLyk} is to perform a Fej\'er-Riesz factorization of a non-negative trigonometric polynomial, which requires an iterative process,  whereas, in contrast, the construction of $x$ in this paper is purely algebraic and can be carried out entirely in rational arithmetic. The algorithm for the solution of the royal tetra-interpolation problem is presented in Section \ref{algorithm}.

The authors are grateful to Nicholas Young for some helpful suggestions.

\section{The phasar derivatives of $\Psi_{\omega} \circ x $ and $\Upsilon_{\omega} \circ x $ } \label{PhasarSection}

The tetrablock was introduced in \cite{AWY}, and it is related to the $\mu_{\mathrm{Diag}}$-synthesis problem. By \cite[Theorem 2.9]{AWY},
	$\mathcal{{E}} \cap \mathbb{R}^3$ is the open tetrahedron with vertices $(1, 1, 1)$, $(1, -1, -1)$, $(-1, 1, -1)$ and $(-1, -1, 1)$.

The following functions are important in the study of $\mathcal{E}$.
\begin{defn} {\em \cite[Definition 2.1]{AWY} } \label{DefOfUps}
For $ x=(x_{1},x_{2},x_{3}) \in \mathbb{C}^3$ and $z \in \mathbb{C}$ we define 
\begin{equation}
\Psi(z,x) =\frac{x_{3}z-x_{1}}{x_{2}z-1} \ \  \textnormal{when} \quad x_{2}z-1 \neq 0,
\end{equation}
\begin{equation} 
\Upsilon(z,x) =\frac{x_{3}z-x_{2}}{x_{1}z-1} \ \ \textnormal{when} \quad x_{1}z-1 \neq 0.
\end{equation} 
{\em For $\omega \in \mathbb{T} $, let}
\begin{equation}
\displaystyle \Psi_{\omega}(x) =\frac{x_{3} \omega - x_{1}}{x_{2} \omega-1} \ \textnormal{when} \  x_{2} \omega-1 \neq 0.
\end{equation}
\end{defn}

Further we will need the following description of the tetrablock from \cite{AWY}.
\begin{thm} {\em \cite[Theorem 2.4]{AWY} } \label{DefOfTetra}
	For $ x \in \mathbb{C}^3 $ the following are equivalent.
	\begin{enumerate}
		\item[{\em (i)}] $ x \in {\bar{\mathcal{E}}}$;
		\item[(ii)] $ ||\Psi( . ,x) ||_{H^\infty} \leq  1 $ and if $x_1 x_2 = x_3 $ then, in addition, $ | x_2| \leq 1$;
		\item[{\em (iii)}] $ ||\Upsilon( . ,x) ||_{H^\infty} \ \leq \ 1 $ \ and if \ $x_1 x_2 = x_3 $ then, in addition, $ | x_1| \ \leq 1$;
		\item[{\em (iv)}] $ |x_1|^{2} - |x_2|^{2} + |x_3|^{2} + 2 | x_2 - {\bar{x_1}} x_3 | \leq 1$ and $ | x_2| \leq 1$;
		\item[{\em (v)}] $ | x_1 - {\bar{x_2}} x_3 | + | x_2 - {\bar{x_1}} x_3 | \leq 1 - | x_3|^{2}$ and if $|x_3| = 1$ then, in addition, $|x_1| \leq 1$.		
	\end{enumerate}	
\end{thm}

By \cite[Theorem 2.9]{AWY},	$\bar{\mathcal{E}}$ is polynomially convex. 
Since $\bar{\mathcal{E}}$ is polynomially convex, there is a smallest closed boundary 
$ b \bar{\mathcal{E}}$ of $ \overline{\mathcal{E}}$, which is called 
{\em  the distinguished boundary of $  \overline{\mathcal{E}}$}. 

\begin{thm} \label{PropOfTetra}
	{\em \cite[Theorem 7.1] {AWY}}  For x $\in$ $\mathbb{C}^3$ the following are equivalent.
	\begin{enumerate}
\item[\em(i)] $ x_1 = \bar{ x}_2x_3, |x_3| = 1$ and  $|x_2| \leq1  $;
		\item[\em (ii)] either $x_1x_2 \neq x_3 $ and $\Psi (.,x) $ is an automorphism of $\mathbb{D}$ or $ x_1x_2=x_3$ and $|x_1| =|x_2| = |x_3| =1$;
		\item[\em (iii)] x is a peak point of $\mathcal{\bar{E}}$;
		\item[\em (iv)] there exists a $2 \times 2$ unitary matrix U such that $x=\pi(U)$ where 
		$$
		\pi : \mathbb{C}^{2 \times 2} \rightarrow \mathbb{C}^{3} : U = [u_{ij}] \mapsto (u_{11}, u_{22}, \det U); 	$$ 
		\item[\em (v)] there exists a symmetric $2 \times 2$ unitary matrix  U such that $x=\pi(U)$;
		\item[\em (vi)] $x \in b\bar{\mathcal{E}}$;
		\item[\em (vii)] $x\in \mathcal{\bar{E}}$ and $|x_3|=1$.
	\end{enumerate}
\end{thm}
\noindent By \cite[Corollary 7.2] {AWY},
 $ b\bar{\mathcal{E}} $ is homeomorphic to $ \overline{ \mathbb{D}} \times \mathbb{T} $.

\begin{lem} \label{InnerUppsi}
	Let $(x_1, x_2, x_3) \in \mathcal{\overline{E}}$ be such that $x_1 x_2 \neq x_3$. For any $\omega \in \mathbb{T}$,
	$$
	|\Psi_{\omega}(x_1, x_2, x_3)| = 1 \ \text{if and only if } \  2 \omega (x_2 - \overline{x_1} x_3 ) = 1  - |x_1|^{2} + |x_2|^{2} - |x_3|^{2}.
	$$
\end{lem}	
	\begin{proof} Consider $\omega \in \mathbb{T}$.
		\begin{align} \nonumber
		|\Psi_{\omega}(x_1, x_2, x_3)| = 1  &\Leftrightarrow \left|\frac {  x_3 \omega - x_1 }  { x_2  \omega - 1  } \right| = 1 \\ \nonumber
	&\Leftrightarrow |\omega x_3 - x_1 |^{2} = |x_2  \omega - 1 |^{2} \\\nonumber
		&\Leftrightarrow |x_3|^{2} - 2 \text{Re}  (\omega x_3 \overline{x_1} ) + |x_1|^{2} = |x_2|^{2} - 2 \text{Re}(x_2 \omega) +1  \\ 
\nonumber
		&\Leftrightarrow |x_1|^{2} - |x_2|^{2} +  |x_3|^{2} + 2\text{Re} ( \omega (x_2 - \overline{x_1} x_3 )) = 1 \\ \label{eqquu}
		&\Leftrightarrow 2 \text{Re}( \omega (x_2 - \overline{x_1} x_3 )) = 1 - |x_1|^{2} + |x_2|^{2} -  |x_3|^{2}.  
		\end{align} 
Since $(x_1, x_2, x_3) \in \mathcal{\overline{E}}$, by \cite[Theorem 2.4 (vii)]{AWY} (see Theorem \ref{DefOfTetra}), 
$$
		2 | x_2 - {\bar{x_1}} x_3 | \leq 1 - |x_1|^{2} + |x_2|^{2} - |x_3|^{2},
		$$ and $|x_2| \leq 1$. Therefore,
		$$
		 2 \text{Re} ( \omega (x_2 - \overline{x_1} x_3 )) \leq 2 | x_2 - \overline{x_1} x_3 | \leq 1 - |x_1|^{2} + |x_2|^{2} - |x_3|^{2} = 2 \text{Re} ( \omega (x_2 - \overline{x_1} x_3 )).
		$$ 
		Thus
		$$
		2 \text{Re} ( \omega (x_2 - \overline{x_1} x_3 )) = 2 | x_2 - \overline{x_1} x_3 | = 1 - |x_1|^{2} + |x_2|^{2} - |x_3|^{2} .
		$$
		Hence, for any $\omega \in \mathbb{T}$, $|\Psi_{\omega}(x_1, x_2, x_3)| = 1$  if and only if 
		$$
		2 \omega (x_2 - \overline{x_1} x_3 ) = 1  - |x_1|^{2} + |x_2|^{2} - |x_3|^{2}.
		$$ 
\end{proof}

For a rational $\overline{\mathcal{{E}}}$-inner function $x= (x_1, x_2,x_3) : \mathbb{D} \rightarrow \overline{\mathcal{E}}$, we consider the rational functions $ \psi_{\omega} : \mathbb{D}  \rightarrow  \overline{\mathbb{D}} $ and $\upsilon_{\omega} : \mathbb{D}  \rightarrow  \overline{\mathbb{D}} $ given, for any $\omega \in \mathbb{T}$, by
$$
\psi_{\omega}(\lambda) = \Psi_{\omega} \circ x (\lambda) \ = \ \frac { \omega x_3 - x_1 }  { x_2  \omega - 1  } (\lambda), \quad \text{for all} \ \lambda \in \overline{\mathbb{D}} \ \ \text{such that} \  x_{2}(\lambda) \omega - 1 \neq 0,
$$ 
and
$$
\upsilon_{\omega}(\lambda) = \Upsilon_{\omega} \circ x (\lambda) =  \frac{x_{3} \omega -x_{2}}{x_{1}\omega-1}(\lambda), \quad  \text{for all} \  \lambda \in \overline{\mathbb{D}} \ \ \text{such that} \ 
x_{1}(\lambda) \omega - 1 \neq 0.
$$  

In  \cite{ALY2} we introduced the terminology of  { \it the phasar derivative} $Af(z)$ for any differentiable function $f : \mathbb{T} \rightarrow \mathbb{C}  \setminus  \{0\} $ at $ z = e^{i\theta} \in \mathbb{T} $ and wrote down some useful elementary properties of phasar derivatives, see Definition \ref{phasar-derivative}.
  
\begin{prop} \label{prop2.4Finite} \cite{ALY2}					
	\begin{enumerate}
		\item[{\em (i)}] For differentiable functions $ \psi, \varphi : \mathbb{T} \rightarrow \mathbb{C} \setminus \{ 0 \}$ and for any $ c \in \mathbb{C} \setminus \{ 0 \}$,
		\begin{equation}
		A(\psi \varphi) = A \psi + A\varphi \quad and \quad A(c\psi) = A\psi.
		\end{equation}							
		\item[{\em (ii)}] For any rational inner function $\varphi$ and for all $ z \in \mathbb{T}$, 
		\begin{equation}
		A\varphi(z) = z \frac{\varphi^{\prime}(z)} {\varphi(z)}.
		\end{equation}							
		\item[{\em (iii)}] If $ \alpha \in \mathbb{D} $ and 
		$$ 
		B_a(z)=  \frac {z-\alpha} {1-\bar{\alpha}z},
		$$
		then
		$$
		AB_{\alpha}(z) = \frac { 1-|\alpha|^{2}} {|z-\alpha|^{2}}\ > \ 0 \quad for \ z \in \mathbb{T}.
		$$
		\item[{\em (iv)}] For any rational inner function $p$, 
		$$
		Ap(z) \ > \ 0 \quad for \ all \ z \in \mathbb{T} .
		$$
	\end{enumerate}
\end{prop} 
Let us recall that, by definition, $\sigma \in \mathbb{T}$ is a royal node of a tetra-inner function $x = (x_1, x_2, x_3)$ if $x_3(\sigma) - x_1(\sigma) x_2(\sigma) = 0$.    

\begin{lem} \label{sigmalemma}
	Let $ x = (x_1, x_2, x_3) $ be a rational tetra-inner function and let $\sigma \in \mathbb{T} $ be a royal node of $x$. Then $\sigma$ is a zero of the function  $ x_3 - x_1 x_2 $ of multiplicity at least $2$. 
\end{lem}		
\begin{proof}
	If $\lambda \in \mathbb{T}$, we have $ x_3(\lambda) - x_1(\lambda) x_2(\lambda) = 0 $ if and only if $\lambda$ is a royal node of $x$. 
	
	For $ \lambda \in \mathbb{T}$, since $x$ is a tetra-inner function, by Theorem \ref{PropOfTetra} (i),  
	\begin{eqnarray}\label{fge0}
		\overline{x_3(\lambda)} (x_3(\lambda) - x_1(\lambda) x_2(\lambda))   &=& 
		\overline{x_3(\lambda)}x_3(\lambda) -  \overline{x_3(\lambda)} x_1(\lambda) x_2(\lambda) \nonumber\\
		&=&  |x_3(\lambda)|^{2} - \overline{x_3(\lambda)} ( \overline{x_2(\lambda)} x_3(\lambda ) ) x_2(\lambda) \nonumber\\
		&=& 1 - |x_3(\lambda)|^{2} |x_2(\lambda)|^{2} \nonumber \\
		&=& 1 - | x_2(\lambda) |^{2} \geq 0, \ \text{ since} \ |x_2(\lambda)| \leq 1  \text{ on} \ \mathbb{T}.
	\end{eqnarray}
For $\theta \in \R$, let  $$	f(\theta) = \overline{x_3(e^{i\theta})} (x_3(e^{i\theta}) - x_1(e^{i\theta}) x_2(e^{i\theta})),$$
and let $ \sigma = e^{i\xi} $.
	By assumption, $ x_3(\sigma) - x_1(\sigma) x_2(\sigma) = 0 $,
and so $f(\xi) = 0$.
By inequality \eqref{fge0}, the function 
	$$
	f(\theta) = \overline{x_3(e^{i\theta})} (x_3(e^{i\theta}) - x_1(e^{i\theta}) x_2(e^{i\theta})) =
	1 - | x_2(e^{i\theta})|^{2} \geq 0, $$
and so it	has a local minimum at $\xi$. Therefore $\xi$ is a critical point of $f$, and $$ \frac{d}{d \theta} ( 1 - | x_2(e^{i\theta})|^{2} )_{|\xi} = 0.$$
Thus,\begin{eqnarray*}
		0 &=& 
\frac{d}{d \theta} ( 1 - | x_2(e^{i\theta})|^{2} )_{|\xi} \\ 
		 &=& 	\frac{d}{d \theta} ( \overline{x_3(e^{i\theta})} (x_3(e^{i\theta}) - x_1(e^{i\theta}) x_2(e^{i\theta}))_{|\xi}\\
		&=& \frac{d}{d \theta} (\overline{x_3(e^{i\theta})})_{|\xi}  (x_3(e^{i\theta}) - x_1(e^{i\theta}) x_2(e^{i\theta}))_{|\xi} + {\overline{x_3(e^{i\theta})}}_{|\xi} \frac{d}{d \theta} (x_3(e^{i\theta}) - x_1(e^{i\theta}) x_2(e^{i\theta}))_{|\xi}  \\
		&=& \frac{d}{d \theta} (\overline{x_3(e^{i\theta})})_{|\xi} \times 0 +\overline{x_3(e^{i\xi})} \bigg[ \frac{d}{d \theta} {x_3(e^{i\theta})}_{|\xi} - \frac{d}{d \theta} ( x_1(e^{i\theta}) x_2(e^{i\theta}) )_{|\xi} \bigg] \\
		&=& \overline{x_3(e^{i\xi})} \bigg[ ie^{i\xi}  x_3'(e^{i\xi}) - (x_1(e^{i\xi}) ie^{i\xi} x_2'(e^{i\xi})+  ie^{i\xi} x_1'(e^{i\xi})  x_2(e^{i\xi}))  \bigg].
	\end{eqnarray*}
	Note that $|x_3(e^{i\xi})| = 1 $, hence $ x_3'(\sigma) = x_1(\sigma) x_2'(\sigma)+ x_1'(\sigma)  x_2(\sigma)$. Thus 
$$ x_3(\sigma) - x_1(\sigma) x_2(\sigma) = 0 \ \text{ and} \  (x_3(\sigma) - x_1(\sigma) x_2(\sigma))' = 0. $$ Therefore $\sigma$ is a zero of $(x_3 - x_1 x_2) $ of multiplicity at least $2$. 
\end{proof}	
\begin{prop} \label{PhasarTheoremForTetra}
		Let $ x = (x_1, x_2, x_3) $ be a rational $\overline{\mathcal{{E}}}$-inner function. Let $\sigma \in \mathbb{T} $ be a royal node of $x$.  Suppose  
$x(\sigma) = (\eta, \tilde{\eta}, \eta \tilde{\eta}), \  \omega \in \mathbb{T} \  \text{and} \  \omega \tilde{\eta} \neq 1 $. Then
	$$
	A(\Psi_{\omega} \circ x) (\sigma) = A x_1(\sigma). 
	$$
\end{prop}
\begin{proof}
	Since $x$ is a rational $\overline{\mathcal{E}}$-inner function, for almost all $\lambda \in \mathbb{T}, \ x(\lambda) \in b\overline{\mathcal{{E}}}$, and, by Theorem \ref{PropOfTetra}, for almost all $\lambda \in \mathbb{T}$,  $ x_1(\lambda) = \overline{ x_2(\lambda)} x_3(\lambda), |x_3(\lambda)| = 1$ and  
$|x_2(\lambda)| \leq   1$.                                                                                                                                                                                                                                                                                                                                                                                                                                                                                                                                                                                                                                                                                                                                                                                                                                                                                                                                                                                                                                                                                                                                                                                                                                                                                                                                                                                                                                                                                                                                                                                                                                                                                                                                                                                                                                                                                                                                                                                                                                                                                                                                                                                                                                                                                                                                                                                                                                                                                                                                                                                                                                                                                                                                                                                                                                                                                                                                                                                                                                                                                                                                                                                                                                                                                                                                                                                                                                                                                                                                                                                                                                                                                                                                                                                                                                                                                                                                                                                                                                                                                                                                                                                                                                                                                                                                                                                                                                                                                                                                                                                                                                                                                                                                                                                                                                                                                                                                                                                                                                                                                                                                                                                                                                                                                                                                                                                                                                                                                                                                                                                                                                                                                                                                                                                                                                                                                                                                                                                                                                                                                                                                                                                                                                                                                                                                                                                                                                                                    By Proposition \ref{prop2.4Finite}, for every $z \in \mathbb{T}$, and every rational inner function $\varphi$,
	$$
	A\varphi(z) = z \frac{\varphi'(z)}{\varphi(z)}.
	$$
	For $\sigma \in \mathbb{T}$ such that $x(\sigma) \in \mathcal{R_{\overline{\mathcal{E}}}}$ and $\omega \tilde{\eta} \neq 1$,	
	\begin{eqnarray*}
		A(\Psi_{\omega} \circ x) (\sigma) &=& A(\omega x_3 - x_1)(\sigma) - A( x_2 \omega - 1 )(\sigma) \\
		&=&\sigma \frac{(\omega x_3 - x_1)' (\sigma)}{(\omega x_3 - x_1)(\sigma)} - \sigma \frac{(x_2 \omega - 1)' (\sigma)}{(x_2 \omega - 1)(\sigma)} \\
		&=&\frac{\sigma}{\omega \tilde{\eta} - 1} \Big(\frac{\omega x_3'(\sigma) - x_1'(\sigma)}{\eta} - \omega x_2'(\sigma)\Big).
	\end{eqnarray*}
	Since $x_3(\sigma) \in \mathcal{R_{\bar{\mathcal{E}}}},$ we have $ x_3(\sigma) = x_1(\sigma) x_2(\sigma)$, and, by Lemma \ref{sigmalemma}, $\sigma$ is a zero of $ x_3 - x_1 x_2 $ of multiplicity at least 2. Thus $(x_3 - x_1 x_2 )' (\sigma) = 0 $ and 
	\begin{equation} \label{royal_relation1}
	x_3'(\sigma) = x_1(\sigma) x_2'(\sigma) + x_2(\sigma) x_1'(\sigma) = \eta x_2'(\sigma) + \tilde{\eta} x_1'(\sigma).
	\end{equation}
Thus, by equation ({\ref{royal_relation1}}), we have	
	\begin{eqnarray*}
		A(\Psi_{\omega} \circ x) (\sigma)	
		&=& \frac{\sigma}{\omega \tilde{\eta} - 1} \Big(\frac{\omega(\eta x_2'(\sigma) + \tilde{\eta} x_1'(\sigma)) - x_1'(\sigma) }{\eta} - \omega x_2'(\sigma)\Big) 	\\	
		&=&\frac{\sigma}{\omega \tilde{\eta} - 1}  \Big(\frac{\omega \eta x_2'(\sigma) + x_1'(\sigma)(\omega \tilde{\eta} - 1 ) - \eta \omega x_2'(\sigma)
		}{\eta } \Big) \\
		&=& \sigma \Big( \frac{x_1'(\sigma)}{\eta} \Big) = \sigma \Big( \frac{x_1'(\sigma)}{x_1(\sigma)} \Big) = A x_1(\sigma).
	\end{eqnarray*}
\end{proof}

\begin{prop} \label{PhasarTheoremForTetra2}	
	Let $ x = (x_1, x_2, x_3) $ be a rational $\overline{\mathcal{{E}}}$-inner function. Let $\sigma \in \mathbb{T} $ be a royal node of $x$.  Suppose  $x(\sigma) = (\eta, \tilde{\eta}, \eta \tilde{\eta}), \  \omega \in \mathbb{T} \  \textnormal{and} \  \omega \eta \neq 1 .$ Then
	$$
	A(\Upsilon_{\omega} \circ x) (\sigma) = A x_2(\sigma). 
	$$
\end{prop}
\begin{proof}
	Since $x$ is a rational $\overline{\mathcal{E}}$-inner function, then for almost all $\lambda \in \mathbb{T}, \ x(\lambda) \in b\overline{\mathcal{{E}}}$, and, by Theorem \ref{PropOfTetra} (i ), for almost all $\lambda \in \mathbb{T}$,  $ x_1 = \bar{ x}_2x_3, |x_3| = 1$ and  $|x_2| \leq1  $.
	By Proposition \ref{prop2.4Finite}, for every $z \in \mathbb{T}$, and every rational inner function $\varphi$,
	$$
	A\varphi(z) = z \frac{\varphi'(z)}{\varphi(z)}.
	$$
	For $\sigma \in \mathbb{T}$ such that $x(\sigma) \in \mathcal{R_{\bar{\mathcal{E}}}}$, and $\omega \eta \neq 1$,	
	\begin{eqnarray*}
		A(\Upsilon_{\omega} \circ x) (\sigma) &=& A(\omega x_3 - x_2)(\sigma) - A(  \omega x_1 - 1 )(\sigma) \\
		&=&\sigma \frac{(\omega x_3 - x_2)' (\sigma)}{(\omega x_3 - x_2)(\sigma)} - \sigma \frac{(x_1 \omega - 1)' (\sigma)}{(x_1 \omega - 1)(\sigma)} \\
		&=&\frac{\sigma}{\omega \eta - 1} \Big(\frac{\omega x_3'(\sigma) - x_2'(\sigma)}{\tilde{\eta}} - \omega x_1'(\sigma)\Big).
	\end{eqnarray*}
	Since $x_3(\sigma) \in \mathcal{R_{\bar{\mathcal{E}}}},$ we have $ x_3(\sigma) = x_1(\sigma) x_2(\sigma)$, and, by Lemma \ref{sigmalemma}, $\sigma$ is a zero of $ x_3 - x_1 x_2 $ of multiplicity at least 2. Thus $(x_3 - x_1 x_2 )' (\sigma) = 0 $ and 
	\begin{equation} \label{royal_relation}
	x_3'(\sigma) = x_1(\sigma) x_2'(\sigma) + x_2(\sigma) x_1'(\sigma) = \eta x_2'(\sigma) + \tilde{\eta} x_1'(\sigma).
	\end{equation}
Thus, by equation ({\ref{royal_relation}}), we have 
\begin{eqnarray*}
		A(\Upsilon_{\omega} \circ x) (\sigma)	
		&=& \frac{\sigma}{\omega \eta - 1} \Big(\frac{\omega(\eta x_2'(\sigma) + \tilde{\eta} x_1'(\sigma)) - x_2'(\sigma) }{\tilde{\eta}} - \omega x_1'(\sigma)\Big) 	\\	
		&=&\frac{\sigma}{\omega \eta - 1}  \Big(\frac{\omega \eta x_2'(\sigma) + x_1'(\sigma) \omega \tilde{\eta} -  x_2'(\sigma) - \tilde{\eta} \omega x_1'(\sigma)
		}{\tilde{\eta} } \Big) \\
	    &=& \frac{\sigma}{\omega \eta - 1} \Big( \frac{x_2'(\sigma) (\omega \eta -1)}{\tilde{\eta}} \Big) = \sigma \Big( \frac{x_2'(\sigma)}{\tilde{\eta}} \Big) \\	
		&=& \sigma \Big( \frac{x_2'(\sigma)}{x_2(\sigma)} \Big) = A x_2(\sigma) .
	\end{eqnarray*}
\end{proof}

\section{Rational tetra-inner functions and royal polynomials } \label{RoyalPolSection}

In this section we will show how to construct rational $\bar{\mathcal{E}}$-inner functions with prescribed royal nodes and values. To describe this construction we need several theorems and definitions from   {\cite {AlsLyk}}. Detailed proofs of these statements are given in {\cite {AlsLyk,OZY}}.
	
For a polynomial $p$ of degree less than or equal to $n$, where $n \ge 0$, we define the polynomial $p^{\sim n }$ by $$p^{\sim n }(\lambda) = \lambda^{n}  \overline{p(\displaystyle{\frac{1}{\overline{\lambda}}})}.$$	

\begin{thm} {\em \cite[Theorem 4.15] {AlsLyk} } \label{rat-tetra-function}
	If $ x = (x_1, x_2, x_3) $ is a rational $ \bar{\mathcal{E}}$-inner function of degree $n$, then there exist polynomials $ E_1, E_2, D $ such that 
	\begin{enumerate}
		\item[\em (i)]  {\em {deg}}($E_1$), {\em {deg}}($E_2$), {\em {deg}}(D) $\leq n $,
		\item[\em (ii)]  $ D(\lambda) \neq \  0 $ on $ \overline{\mathbb{D}}$,
		\item[\em (iii)] $ E_1(\lambda) = E_2^{\sim n }(\lambda) , \ \mbox{for all} \  \lambda \in \mathbb{T}$, 
		\item[\em (iv)] $ | E_i(\lambda)| \leq | D(\lambda)| $ on $ \overline{\mathbb{D}} , \ i = 1,2$,
		\item[\em (v)] $ x_1 = \frac{E_1} {D}$ on $ \overline{\mathbb{D}}$,
		\item[\em (vi)] $ x_2 = \frac{E_2} {D}$ on $  \overline{\mathbb{D}}$,
		\item[\em (vii)] $ x_3 = \frac{D^{\sim n }} {D}$ on $  \overline{\mathbb{D}}$.	
	\end{enumerate}
\end{thm}

	\begin{defn}\label{RoyalPolyDef}
	Let $ x = (x_1, x_2, x_3)$ be a rational tetra-inner function of degree $n$. The {\em royal polynomial } of $x$ is 
	$$
	R_{x}(\lambda) = D(\lambda) D^{\sim n }(\lambda) - E_1(\lambda) E_2(\lambda),
	$$ where $ E_1, E_2, D $ are as in {\em Theorem  \ref{rat-tetra-function}}.
	\end{defn}
	\begin{rem} \label{DefOfRoyalNodes}
	For a rational tetra-inner function $x$, since $ D(\lambda) \neq   0 $ on $ \overline{\mathbb{D}}$, the zeroes of $R_{x}$ are the zeroes of the function $x_3 - x_1 x_2$. 
\end{rem}			
\begin{lem} \label{LemmaForUnitDisc}
		Let  $ x = (x_1, x_2, x_3)$ be a rational $\overline{\mathcal{E}}$-inner function, and let $\sigma \in \bar{\mathbb{D}}$ be a royal node of $x$. 
		If $\sigma \in \mathbb{T}$, then $|x_1(\sigma)| = 1$ and $|x_2(\sigma)| = 1$.
	\end{lem}
\begin{proof}
Since $x$ is an $\overline{\mathcal{E}}$-inner function, by definition  of  tetra-inner functions, $x(\sigma) \in b \overline{\mathcal{{E}}}$ \ for $\sigma \in \mathbb{T}$. By Theorem \ref{PropOfTetra}, $x_1(\sigma) = \overline{x_2(\sigma)} x_3(\sigma)$,  $|x_3(\sigma)| = 1$ and $|x_2(\sigma)| \leq 1$. By assumption $\sigma$ is a royal node of $x$. Thus $ x_3(\sigma) = x_1(\sigma) x_2(\sigma) $, and so   $|x_1(\sigma)| = 1$ and $|x_2(\sigma)| = 1$ since  $|x_3(\sigma)| = 1$.  
\end{proof}

\begin{definition}  \cite[Definition 3.4]{ALY4} \label{n-sym}
	We say that a polynomial $f$ is {\em $n$-symmetric}  if $\deg(f) \leq n$ and $f^{\sim n}=f$. 
\end{definition}

\begin{prop} {\em \cite {AlsLyk}} 
	Let $x$ be a rational $\overline{\mathcal{E}}$-inner function of degree $n$ and let $R_{x}$ be the royal polynomial of $x$. Then $	R_{x}$ is $2$n-symmetric and the zeros of $	R_{x}$ on $\mathbb{T}$ have even order or infinite order.    
\end{prop}
	\begin{defn} \label{DefOfn,kType}
	Let  $ x = (x_1, x_2, x_3)$ be a rational $\bar{\mathcal{E}}$-inner function such that $x(\overline{\mathbb{D}} )  \nsubseteq \mathcal{R}_{\overline{\mathcal{E}}} $ and let $R_x$ be  the royal polynomial of $x$. If $\sigma$ is a zero of $R_{x}$ of order $\ell$, we define \index{multiplicity of a royal node } the multiplicity $\# \sigma$ of $\sigma $ (as a royal node of $x$) by
	$$
    \#\sigma \ \ \ = \ \ \ \begin{cases}
				\ell & \text{if } \ \sigma \in \mathbb{D}\\
				\frac {1}{2}\ell & \text{if} \ \ \sigma \in \mathbb{T}.
				\end{cases}
				$$
				We define the {\em type} of $x$ to be the ordered pair $(n,k)$, where $n$ is the number of royal nodes of $x$ that lie in $\bar{\mathbb{D}},$ counted with multplicity, and $k$ is the number of royal nodes of $x$ that lie in $\mathbb{T}$, counted with multiplicity. 
		 $\mathcal{R}^{n,k}$  denotes the collection of rational $\bar{\mathcal{E}}$-inner functions of type $(n,k)$.
			\end{defn}	
		\begin{defn}{\em \cite {AlsLyk}}
			The {\em  degree} of a rational $\overline{\mathcal{E}}$-inner function $x$, denoted by {\em{deg}}$(x)$ is defined to be $x_*(1)$, where $x_* : \mathbb{Z} = \pi_1(\mathbb{T}) \rightarrow \pi_1(b\overline{\mathcal{E}}) $ is the homomorphism of fundamental groups induced by $x$ when $x$ is regarded as a continuous map from $\mathbb{T} $ to $b\overline{\mathcal{E}}$. 
		\end{defn}
	\begin{prop}{\em \cite {AlsLyk}}
		For any rational $\overline{\mathcal{E}}$-inner function $x$, {\em{deg}}$(x)$ is the degree {\em{deg}}$(x_3)$ (in the usual sense) of the finite Blaschke product $x_3$.
	\end{prop}
		\begin{thm} {\em \cite {AlsLyk}} \label{Degreeofx}
			If $ x \in \mathcal{R}^{n,k} $ is non-constant, then the degree of $x$ is equal to $n$. 
		\end{thm}
	
		\begin{thm} \cite {AlsLyk}
		Let $ x$ be a non-constant rational $\overline{\mathcal{{E}}}$-inner function of degree $n$. Then, either $x(\overline{\mathbb{D}}) \subseteq \mathcal{R}_{\overline{\mathcal{E}}} $ or $x(\overline{\mathbb{D}})$ meets $\mathcal{R}_{\overline{\mathcal{E}}}$ exactly $n$ times.
	\end{thm}
\begin{prop} \label{DegThm}
	Let $x =(x_1, x_2, x_3)$ be a non-constant rational $\overline{\mathcal{{E}}}$-inner function and let $\omega \in \mathbb{T}$ be such that $\omega x_2(\lambda) - 1 \neq 0$ for all $\lambda \in \mathbb{D}$. Then the rational function $\Psi_{\omega} \circ x = \displaystyle{\frac{\omega x_3 - x_1}{x_2 \omega - 1}}$ has a cancellation at $\zeta \in \overline{\mathbb{D}}$ if and only if the following conditions are satisfied : $\zeta \in \mathbb{T}, \  \zeta $ is a royal node of $x$ and $\omega = \overline{x_2(\zeta)} $. 
\end{prop}
\begin{proof}
	 Let $\zeta \in \mathbb{T}$ be a royal node of $x$ such that $x(\zeta) = (\eta, \tilde{\eta} , \eta \tilde{\eta})$. By Lemma \ref{LemmaForUnitDisc}, $|\eta| =1 $ and $|\tilde{\eta}| = 1$. If $\omega= \overline{\tilde{\eta}} \in \mathbb{T}$, then
	$$
	\omega x_3(\zeta) - x_1(\zeta) = \overline{\tilde{\eta}} \eta \tilde{\eta} - \eta = |\tilde{\eta}|^{2} \eta - \eta = \eta - \eta = 0 ,
	$$
	and
	$$
	x_2(\zeta) \omega -1 = \tilde{\eta} \overline{\tilde{\eta}} - 1 = |\tilde{\eta}|^{2} - 1 = 0.
	$$
	Thus, $\Psi_{\omega} \circ x =  \displaystyle{\frac{\omega x_3(\lambda) - x_1(\lambda)}{x_2(\lambda) \omega - 1} }$ has at least one cancellation at such  $\zeta \in \mathbb{T}$.
	
	Conversely, by assumption $\Psi_{\omega} \circ x$ has a cancellation at $\zeta \in \overline{\mathbb{D}}$, and so
	$$
	(\omega x_3 - x_1)(\zeta) = 0 = (x_2 \omega -1)(\zeta). 
	$$
	Therefore, $x_2(\zeta) \omega = 1$ \ and \ $\omega x_3(\zeta) = x_1(\zeta)$. Since $x_2(\zeta) \omega = 1,$  it implies that  $ x_2(\zeta) = \overline{\omega} \in \mathbb{T}$, so $|x_2(\zeta)| = 1$. Since $x_2 : \mathbb{D} \rightarrow \overline{\mathbb{D}}$ is a rational analytic function with  $|x_2(\zeta)| = 1$, by the maximum principle, $\zeta \in \mathbb{T}$, or $x_2(\lambda) = \overline{\omega}$ for all $\lambda \in \overline{\mathbb{D}}$. By  assumption,  $\omega x_2(\lambda) - 1 \neq 0$ for all $\lambda \in \mathbb{D}$. Hence the function $x_2 \neq \overline{\omega}$ on $\mathbb{D}$.  Therefore,  $ \zeta \in \mathbb{T}$. 
	Note
	\begin{align*}
	\omega x_3(\zeta) = x_1(\zeta) \ &\implies \ \overline{x_2(\zeta)}  x_3(\zeta) = x_1(\zeta) \\
	&\implies \ x_3(\zeta) = x_1(\zeta) x_2(\zeta). 
	\end{align*}
	Thus, $\zeta \in \mathbb{T}$ is a royal node for $x$, and  $\omega = \overline{x_2(\zeta)} $.
\end{proof}

\section{Criteria for the solvability of the Blaschke interpolation problem}\label{Criteria for the solvability of the Blaschke}

In this section, for the convenience of the reader, we collect some known facts about finite Blaschke products that we need.  They may be found in several places, but the most economical source for our purposes is \cite{ALY3}, which assembles precisely the results which we require.

As mentioned in the Introduction, there is an extensive literature on boundary interpolation problems. A very valuable source of information about all manner of complex 
interpolation problems is the book of Ball, Gohberg and Rodman \cite{bgr}.
The authors of \cite{BH,Sar,BolKh08,bgr,BD}   make use of Krein spaces, moment theory, measure theory, reproducing kernel theory, realization theory and de Branges space theory.  They obtain far-reaching results, including generalizations to matrix-valued functions and to functions allowed to have a limited number of poles in a disc or half-plane. See also papers of \cite{Bol10,Mar,power} for  elementary treatments of interpolation problems.
The monograph \cite{BD} by Bolotnikov and Dym is entirely devoted to 
boundary  interpolation problems for the Schur class. They reformulate the problem within the framework of the Ukrainian school's Abstract Interpolation Problem
and solve it by means of operator theory in de Branges-Rovnyak spaces.

The Blaschke interpolation Problem \ref{The Blaschke interpolation problem}  as described in \cite {ALY3} is an algebraic variant of the classical Pick interpolation problem. 
 One looks for a Blaschke product of degree $ n $ satisfying $ n $ interpolation conditions, rather than a Schur-class function as in the original Pick interpolation problem. We admit interpolation nodes both in the open unit disc and on the unit circle. There is a criterion for the solvability of the Blaschke interpolation problem in terms of the positivity of a ``Pick matrix" formed from the interpolation data.  To obtain a well-posed problem one  imposes additional interpolation conditions, on phasar derivatives at the interpolation nodes on the unit circle. These bounds on the  phasar derivatives appear as the diagonal entries of the Pick matrix.

\begin{defn}
The Schur class is the set of analytic functions $ S $ from $ \mathbb{D} \ to \ \bar{\mathbb{D}}$, $ S : \mathbb{D} \rightarrow
 \bar{\mathbb{D}}  $.
\end{defn}
\begin{defn}
	A function $ f: \mathbb{D} \rightarrow \overline{\mathbb{D}}$ is {\em inner} if it  is an analytic map such that the radial limit
	$$
	\lim_{r \rightarrow 1^{-}} f(r\lambda)  \ \mbox{exists and belongs to }  \mathbb{T}
	$$  
	for almost all $\lambda \in \mathbb{T}$ with respect to Lebesgue measure.
\end{defn}

\begin{defn}
The {\em Pick matrix} associated with Blaschke interpolation data $(\sigma , \eta , \rho ) $  of type $(n,k)$ is the $ n \times n $ matrix $ M = [ m_{ij}]_{i,j=1}^n $ where
$$
m_{ij}=\begin{cases}
    \rho_i & \text{if $i = j \leq k$}.\\
    \displaystyle{\frac {1-\overline{\eta_i}\eta_j}{1-\overline{\sigma_i}\sigma_j}} & \mathrm{otherwise}.
  \end{cases}
$$
\end{defn}

\begin{defn}
	The Pick matrix $ M = [ m_{ij}]_{i,j=1}^n $ is {\em minimally positive} if $M \geq 0  $ and there is no strictly positive diagonal $n \times n$ matrix $D$ such that $M \geq D$. 
	
\end{defn}

The following is a refinement of the Sarason Interpolation Theorem \cite{Sar}.
\begin{thm}  \cite[Theorem 3.3] {ALY3}
Let $ M $ be the Pick matrix associated with Blaschke interpolation data $( \sigma, \eta, \rho)$ of type $(n,k)$.

\begin{enumerate}
\item[{\em (i)}] There exists a function $ \varphi$ in the Schur class such that
\begin{equation}
\varphi(\sigma_j) = \eta_j \quad for \ j = 1, ... ,n, \label{varphi function1}
\end{equation}
and the phasar derivative $A\varphi(\sigma_j)$ exists and satisfies 
\begin{equation}
A\varphi(\sigma_j) \leq \rho_j  \quad for \ j = 1, ... ,k,\label{varphi function2}
\end{equation}
if and only if $ M \geqslant 0 $;
\item[{\em (ii)}] if $M$ is positive semi-definite and of rank $ r < n $ then there is a unique function $ \varphi $ in the Schur class satisfying conditions {\em (\ref{varphi function1})} and {\em (\ref{varphi function2})} above, and this function is a Blaschke product of degree $ r$;
\item[{\em(iii)}] the unique function $\varphi$ in statement {\em (ii)} satisfies 
\begin{equation}
A\varphi(\sigma_j) = \rho_j  \quad for \ j = 1, ... ,k,
\end{equation}
 if and only if $ M $ is minimally positive.

\end{enumerate}
\end{thm}
In \cite{ALY3} the authors described a strategy for the construction of the general solution of the Blaschke interpolation problem (Problem \ref{The Blaschke interpolation problem}). It is to adjoin an additional boundary interpolation condition $\varphi(\tau) = \zeta$ where $\tau \in \mathbb{T}\ \backslash \ 
\{ \sigma_1, . . ., \sigma_k \} $ and $\zeta \in \mathbb{T} $. This augmented problem has a unique solution. All the solutions of Problem \ref{The Blaschke interpolation problem} are then obtained in terms of a unimodular parameter.      
\begin{lem} {\em \cite [Lemma $3.4$] { ALY3}}\label{prop3.4}
If $ C $ is an $ n\times n $ positive definite matrix, $ u $ is an $ n \times 1 $ column, $ \rho = \langle C^{-1} u, u \rangle $ and the $ ( n + 1 ) \times ( n + 1 ) $ matrix $B$ is defined by 
$$
B = \begin{bmatrix}
    C & u  \\
    u^{*} & \rho \\
\end{bmatrix},
$$
then $B$ is positive semi-definite, rank($B$) = $ n $ and
$$
B  \begin{bmatrix}
    -C^{-1} u  \\
    1 \\
\end{bmatrix} = 0.
$$
\end{lem}

The Pick matrix $ B_{\zeta, \tau} $ of the augmented problem is the $ ( n + 1 ) \times ( n + 1 ) $ matrix, 
\begin{equation}
B_{\zeta, \tau} = \begin{bmatrix}
    M & u_{\zeta, \tau} \\
    u^{*}_{\zeta, \tau} & \rho_{\zeta, \tau} \\
\end{bmatrix}, 
\end{equation} where $$
\rho_{\zeta, \tau} = \langle M^{-1} u_{\zeta, \tau}, u_{\zeta, \tau} \rangle,
$$
$ M $ is the Pick matrix associated with  Problem \ref{The Blaschke interpolation problem}, and $u_{\zeta, \tau} $ is the $ n \times 1 $ column matrix defined by 
\begin{equation} \label{ equ lem 3.4}
u_{\zeta, \tau} = \begin{bmatrix}
    \frac { 1-\bar{\eta}_{1}\zeta} { 1-\bar{\sigma}_{1}\tau}  \\ \vdots \\ 
    \frac { 1-\bar{\eta}_{n}\zeta } { 1-\bar{\sigma}_{n}\tau } \\
\end{bmatrix}.
\end{equation}

\begin{thm} {\em \cite [Proposition 3.6] {ALY3}}\label{prop3.6}
If the Pick matrix $M$ associated with {\em Problem \ref{The Blaschke interpolation problem} } is positive definite then, for any $\tau \in \mathbb{T} \ \backslash \ 
 \{ \sigma_1, . . ., \sigma_k \} $ and $\zeta \in \mathbb{T} $, there is at most one solution $\varphi$ of {\em Problem \ref{The Blaschke interpolation problem} }  for which  $\varphi(\tau) = \zeta$.
\end{thm}

The $j$th standard basis vector in $\mathbb{C}^{n}$ will be denoted by $e_j$.

\begin{thm}  \cite [Proposition 3.7] {ALY3}\label{prop3.7}
If the Pick matrix $M$ for  Problem \ref{The Blaschke interpolation problem} is positive definite, if  $\tau \in \mathbb{T} \setminus
 \{ \sigma_1, . . ., \sigma_k \} $ and $\zeta \in \mathbb{T} $ and  if
\begin{equation}
\langle M^{-1}u_{\zeta,\tau} , e_j \rangle \ \neq \  0 
\end{equation}
for $j = 1, . . . , k$, then there is a unique solution $\varphi$ of  Problem \ref{The Blaschke interpolation problem} that satisfies $ \varphi (\tau) = \zeta$.
\end{thm}

The {\em exceptional set} $ Z_\tau$ for  Problem \ref{The Blaschke interpolation problem}  is defined to be
\begin{equation} \label{exceptionalset}
Z_\tau = \{ \zeta \in \mathbb{T} : \text{ for some } j \mbox{ such that } 1\leq j \leq k,  \langle M^{-1}u_{\zeta,\tau} , e_j \rangle =   0   \}
\end{equation}
Define $ n\times 1 $ vectors $ x_\lambda$ and $y_\lambda $ for $ \lambda \in \overline{\mathbb{D}} \backslash \{ \sigma_1, . . ., \sigma_k \} $ by 
\begin{equation} \label{formula_in3.7}
x_\lambda = \begin{bmatrix}
    \frac { 1} { 1-\bar{\sigma}_{1}\lambda}  \\ \vdots \\ 
    \frac { 1} { 1-\bar{\sigma}_{n}\lambda } \\
\end{bmatrix}, \quad \ \ \ y_\lambda = \begin{bmatrix}
    \frac { \bar{\eta}_{1}} { 1-\bar{\sigma}_{1}\lambda }  \\ \vdots \\ 
    \frac { \bar{\eta}_{n} } { 1-\bar{\sigma}_{n}\lambda  } \\
\end{bmatrix}, 
\end{equation}
so that 
\begin{equation}
u_{\zeta,\tau} = x_\tau - \zeta y_\tau
\end{equation}

\begin{thm}  \cite [Proposition 3.8] {ALY3}\label{prop3.8}
\begin{enumerate}
\item[{\em (i)}] For any $\tau \in \mathbb{T}  \setminus
 \{ \sigma_1, . . ., \sigma_k \} $, if 
$$
\langle x_\tau , M^{-1} e_j \rangle =   0 = \langle y_\tau , M^{-1} e_j \rangle \text{ for some } \ j \mbox{ such that }1\leq j \leq k, 
$$
then $ Z_\tau = \mathbb{T}$. 
\item[{\em (ii)}] There exist uncountably many $ \tau \in \mathbb{T} \setminus
 \{ \sigma_1, . . ., \sigma_k \} $ such that 
$$
\langle x_\tau , M^{-1} e_j \rangle =   0 = \langle y_\tau , M^{-1} e_j \rangle 
$$
does not hold for any j such that  $ 1 \leq j \leq k $. Moreover, for such $\tau$, the set $Z_\tau$ consists of at most $k$ points. 
\end{enumerate}
\end{thm}

\begin{thm} \cite[Theorem $3.9$]{ALY3} \label{prop3.9}
Let the Pick matrix $M$ of  Problem \ref{The Blaschke interpolation problem}  be positive definite, and let $ \tau \in \mathbb{T} \setminus \{ \sigma_1, . . . , \sigma_k \} $ be such that the set
$$
Z_{\tau} = \{ \zeta \in \mathbb{T} : u_{\zeta, \tau } \bot M^{-1} e_j \text{ for some } j \mbox{ such that } 1 \leq j\leq k \}
$$
contains at most $k$ points, where $ u_{\zeta, \tau }$ is defined by equation \eqref{ equ lem 3.4}.
\begin{enumerate}

\item[{\em (i)}] If $ \zeta \in \mathbb{T} \setminus Z_{\tau}$, then there is a unique solution $ \varphi_\zeta $ of  Problem \ref{The Blaschke interpolation problem}  that satisfies $ \varphi_\zeta(\tau) = \zeta$.

\item[{\em (ii)}] There exist unique polynomials $ a_\tau, b_\tau, c_\tau, \ and \ d_\tau$ of degree at most $ n$ such that
\begin{equation} \label{equ_3.11}
\begin{bmatrix}
    a_{\tau}(\tau)       & b_{\tau}(\tau) \\
    c_{\tau}(\tau)       & d_{\tau}(\tau) \\
\end{bmatrix}
=
\begin{bmatrix}
    1 & 0  \\
    0 & 1 \\
\end{bmatrix}
\end{equation}
and, for all $ \zeta \in \mathbb{T}$, if $ \varphi $ is a solution of a {\em Problem \ref{The Blaschke interpolation problem} } such that $ \varphi(\tau) = \zeta$, then 
\begin{equation} \label{equ_3.12}
\varphi(\lambda) = \frac {a_{\tau}(\lambda)\zeta + b_{\tau}(\lambda)}{c_{\tau}(\lambda)\zeta + d_{\tau}(\lambda)}
\end{equation}
for all $ \lambda \in \mathbb{D}$.

\item[{\em (iii)}] If \~{a}, \~{b}, \~{c}, \~{d} are rational functions satisfying the equation 
\begin{equation} \label{equation3.13}
\begin{bmatrix}
    \text{\~{a}}(\tau)     & \text{\~{b}}(\tau)  \\
   \text{\~{c}}(\tau)       & \text{\~{d}}(\tau)  \\
\end{bmatrix}
=
\begin{bmatrix}
    1 & 0  \\
    0 & 1 \\
\end{bmatrix} 
\end{equation}
and such that for three distinct points $ \zeta $ in $ \mathbb{T} \setminus Z_{\tau}$, the equation 
\begin{equation} \label{equation3.14}
\frac {a_{\tau}(\lambda)\zeta + b_{\tau}(\lambda)}{c_{\tau}(\lambda)\zeta + d_{\tau}(\lambda)} = \frac { \text{\~{a}}(\lambda)\zeta + \text{\~{b}}(\lambda)  } {\text{\~{c}}(\lambda)\zeta + \text{\~{d}}(\lambda) }
\end{equation}
holds for all $ \lambda \in \mathbb{D} $, then there exists a rational function $ X $ such that \~{a} = $X$$a_{\tau}$, \~{b} = $X b_{\tau} $, \~{c} = $Xc_{\tau}$ and \~{d} = $ X d_{\tau}$. 
\end{enumerate}
\end{thm}

\begin{remark}\label{explicit-solution} {\rm Let $(\si,\eta,\rho)$ be  Blaschke interpolation data of type $(n,k)$. Suppose the Pick matrix $M$ of this problem is positive definite.
The proof of \cite[Theorem 3.9]{ALY3} gives an {\em explicit}  linear fractional parametrization of the solutions of Problem \ref{The Blaschke interpolation problem}.  As in Theorem \ref{prop3.9} choose 
$\tau \in \T\setminus\{\si_1,\dots,\si_k\}$ such that the set
$Z_\tau$ contains at most $k$ points. Then a normalized linear fractional parametrization of the solution set of Problem \ref{The Blaschke interpolation problem} is
\[
\ph=\frac{a_\tau\zeta+b_\tau}{c_\tau \zeta+d_\tau},
\]
where the polynomials
$ a_\tau$, $b_\tau$, $ c_\tau$ and $ d_\tau$ are defined by the equations 
\begin{equation} \label{abcd}
a_\tau=\pi A, \;b_\tau=\pi B,\; c_\tau=\pi C, \; \text{and} \; d_\tau=\pi D.
\end{equation}
Here (see \cite[Theorem 3.9]{ALY3}),
$$
\pi(\la) = (1-\overline{\tau}\lambda)\prod_{j=1}^n \frac{1-\overline{\sigma_j}\lambda}{1-\overline{\sigma_j}\tau},
$$
and 
\begin{align}
A(\lambda) &= -\langle x_\lambda,M^{-1} x_\tau \rangle +\frac{1}{1-\overline{\tau}\lambda} \label{eq6.90},\\
B(\lambda) &= \langle x_\lambda,M^{-1} y_\tau \rangle \label{eq6.100},\\
C(\lambda) &= -\langle y_\lambda,M^{-1} x_\tau \rangle \label{eq6.110},
\end{align}
 and
\begin{align}
D(\lambda) &= \langle y_\lambda,M^{-1} y_\tau \rangle +\frac{1}{1-\overline{\tau}\lambda}. \label{eq6.120}
\end{align}
Different choices of $\tau$  yield different normalized parametrizations  of the solutions of Problem \ref{The Blaschke interpolation problem}.
}
\end{remark}

\cite[Theorem 3.9]{ALY3} tells us  the following.
\begin{cor}  \cite[Corollary 3.12] {ALY3} \label{CorOfDef10}
Let $  (\sigma, \eta, \rho) $ be Blaschke interpolation data of type $(n,k)$. Suppose the Pick matrix $M$ of this problem is positive definite. There exists a normalized linear fractional parameterization 
$$
\varphi = \frac { a\zeta + b} { c\zeta + d}
$$ 
of the solutions of  Problem \ref{The Blaschke interpolation problem}. Moreover 
\begin{enumerate}
\item[{\em(i)}] at least one of the polynomials $ a, b, c, d$ has degree $n$, 
\item[{\em (ii)}] the polynomials $ a, b, c, d$ have no common zero in $ \mathbb{C}$;
\item[{\em (iii)}] $ |c| \leq |d| $ on $ \overline{\mathbb{D}}$.

\end{enumerate}

\end{cor}

\section{From the royal tetra-interpolation problem to the Blaschke interpolation problem}	\label{FromRoyaltoBlash}

		In this section we show that for Blaschke interpolation data $(\sigma, \eta, \rho)$ of type $(n,k)$ knowledge of a solution $x$ of the royal tetra-interpolation problem $(\sigma, \eta, \tilde{\eta},  \rho)$ for some $\tilde{\eta_{j}} \in \overline{\mathbb{D}}$ allows us to construct a solution of the Blaschke interpolation problem.   
	\begin{thm} \label{4.4Tetra}
	Let $x = (x_1, x_2, x_3) $ be a rational $\overline{\mathcal{E}}$-inner function of type $(n,k)$ having distinct royal nodes $\sigma_1, \sigma_2, ... , \sigma_n$ where $\sigma_1, \sigma_2, ... , \sigma_k \in \mathbb{T} $ and $\sigma_{k+1},..., \sigma_n  \in \mathbb{D}$ and corresponding royal values $ \eta_1, .. , \eta_n$ and $ \tilde {\eta}_1,...,\tilde {\eta}_n$, that is,  $ x(\sigma_j) =(\eta_j,{\tilde{\eta}}_j, \eta_j {\tilde{\eta}}_j)$. Let $\rho_j =  A x_1  (\sigma_j)$ for $ j = 1, 2 , .. , k .$ 	
	\begin{enumerate}
	\item[\em (1)] There exists a rational inner
	 function $\varphi$ that solves the Blaschke interpolation Problem \ref{The Blaschke interpolation problem} for $ (\sigma, \eta, \rho)$, that is, such that $ deg(\varphi) = n $,
	\begin{equation}
	\varphi(\sigma_j) = \eta_j \ \ \text{for} \ j = 1, ..., n
	\end{equation}
	and
	\begin{equation}
	A \varphi(\sigma_j) = \rho_j \ \ \text{for} \ j = 1, ..., k.
	\end{equation}
	Any such function $\varphi$ is expressible in the form $\varphi = \Psi_{\omega} \circ x $ for some $\omega \in \mathbb{T}$.							
	\item[\em (2)] There exist polynomials a, b, c, d of degree at most $n$ such that a normalized parametrization of the solutions of {\em Problem \ref{The Blaschke interpolation problem} } is 
	$$
	\varphi = \frac{ a\zeta + b} { c\zeta + d} , \ \ \zeta \in \mathbb{T} . 
	$$
	\item[\em (3)] For any polynomials a, b, c, d as in $(2)$, there exist  $  x_1^{\circ}, x_2^{\circ}, x_3^{\circ} \in \mathbb{C} $ such that
	\begin{equation} \label{11}
	| x_3^{\circ} | = 1, \ \ \ | x_1^{\circ} | < 1, \ \ \  | x_2^{\circ} | < 1,
	\end{equation}
	\begin{equation} \label{22}
	x_1^{\circ} = \overline{x_2^{\circ}} x_3^{\circ}, 
	\end{equation}
	and moreover, 
	\begin{equation} \label{x_1Equation}
	x_1 = \frac { x_1^{\circ} a + b } { x_1^{\circ} c + d }
	\end{equation}
	\begin{equation} \label{x_2Equation}
	x_2 = \frac{x_3^{\circ} c + x_2^{\circ} d } { x_1^{\circ} c + d }
	\end{equation}
	\begin{equation} \label{x_3Equation}
	x_3 = \frac { x_2^{\circ} b + x_3^{\circ} a} { x_1^{\circ} c + d }.
	\end{equation}							
	\end{enumerate}					
	\end{thm}									
	\begin{proof}					
	(1) For $\omega \in \mathbb{T} $ and for a given  rational $\overline{\mathcal{{E}}}$-inner function $x= (x_1, x_2,x_3) : \mathbb{D} \rightarrow \overline{\mathcal{E}}$, we consider the rational function $ \psi_{\omega} : \mathbb{D}  \rightarrow  \overline{\mathbb{D}} $ given by 
	\begin{equation} \label{key equation}
	\psi_{\omega}(\lambda) = \Psi_{\omega} \circ x (\lambda) \ = \ \frac {  x_3 \omega - x_1 }  { x_2  \omega - 1  } (\lambda).
	\end{equation}							
	Then, if $ \omega \ \neq \ \overline{{\tilde {\eta}}_1} , .... , \overline{{\tilde {\eta}}_k}$,
	\begin{equation} \label{equation11}
	\psi_{\omega} (\sigma_j) \ = \ \frac {  x_3(\sigma_j) \omega - x_1(\sigma_j) }  { x_2(\sigma_j)  \omega - 1  }  \ = \ \frac {  \eta_j {\tilde {\eta}}_j \omega - \eta_j} {{\tilde {\eta}}_j \omega - 1 } \ = \ \eta_j \frac { \omega {\tilde {\eta}}_j - 1  } {{\tilde {\eta}}_j \omega- 1 }  = \eta_j \ \text{for} \ j = 1, ... , n. 
	\end{equation}							
	We claim that, for $ \omega \in \mathbb{T} \setminus \{ \overline{{\tilde {\eta}}_1} , .... , \overline{{\tilde {\eta}}_k}   \}$, the function $ \varphi = \psi_{\omega} $ is a solution of  Problem \ref{The Blaschke interpolation problem}. Let us check that $\varphi$ is an inner function from $ \mathbb{D} $ to $\overline{\mathbb{D}} $. For any $\lambda \in \mathbb{T}$, 
	$$
	\varphi(\lambda) = \psi_{\omega}(\lambda) = \frac { \omega x_3(\lambda) - x_1(\lambda) }  { x_2(\lambda)  \omega - 1  }.
	$$ 
	Since $x$ is a rational $\overline{\mathcal{{E}}}$-inner function, $x(\lambda) \in b\overline{\mathcal{{E}}}$ for almost all $\lambda \in \mathbb{T}$, and, by Theorem \ref{PropOfTetra},  $ x_1(\lambda) = \overline{x_2(\lambda)} x_3(\lambda)$ and $|x_3(\lambda)|= 1$ \ for almost all $ \lambda \in \mathbb{T}$. 
	Thus, for almost all $ \lambda \in \mathbb{T} $,
	$$
		\varphi(\lambda) = \psi_{\omega}(\lambda) = \frac { \omega x_3(\lambda) - \overline{x_2(\lambda)} x_3(\lambda) }  { x_2(\lambda)  \omega - 1  } = \frac{x_3(\lambda) (\omega - \overline{x_2(\lambda))}} {\omega x_2 (\lambda) - 1}.
	$$	
	Hence, for $\lambda \in \mathbb{T} ,$
	$$
	|\varphi(\lambda)| = |x_3(\lambda)| \ \Bigg| \frac{\omega - \overline{x_2(\lambda)}}{\omega x_2 (\lambda) - 1} \Bigg|.
	$$		
	Since $ |x_3(\lambda)| = 1, |\omega| = 1$ and  $  
	| \omega - \overline{x_2(\lambda)} 	| = |\overline{\omega} - x_2(\lambda)|$, we have, for almost all $\lambda \in \mathbb{T}$, 
	
	$$
	\Bigg| \frac{\omega (\overline{\omega} - x_2(\lambda))}{\omega x_2 (\lambda) - 1} \Bigg| = 	\Bigg| \frac{1 - x_2(\lambda) \omega}{-(1 - x_2(\lambda) \omega)} \Bigg| = 1.
	$$
	Therefore, for almost all $\lambda \in \mathbb{T}, |\varphi(\lambda)| = 1$. Hence $\varphi$ is a rational inner function. 
	 
	By equation (\ref{equation11}), $\psi_{\omega} $ takes the required values at $\sigma_{1}, ... , \sigma_{n}.$ By Proposition \ref{PhasarTheoremForTetra},
	
\begin{equation}
A(\Psi_{\omega} \circ x) (\sigma_j) = A x_1(\sigma_j) = \rho_j \ \ \ \  \text{for} \ j = 1, 2, ... ,k.
\end{equation}
Furthermore, deg$(\psi_{\omega})=n$ for $\omega \neq \overline{{\tilde {\eta}}_1} , .... , \overline{{\tilde {\eta}}_k} $. By Theorem \ref{rat-tetra-function}, for a rational $\overline{\mathcal{{E}}}$-inner function $x=(x_1, x_2, x_3)$ such that deg$(x_3)=n$ and if $D$ is the denominator when $x_3$ is written in its lowest terms then $x_1$ and $x_2$ can also be represented as rational functions in which the denominator is $D$. Therefore
\begin{equation} \label{Cancellation}
\text{deg}(\psi_{\omega}) = \text{deg}(x_3) - \# \{ \text{cancellations between} \ \omega x_3 - x_1 \ \text{and} \ x_2 \omega -1 \}.
\end{equation}
By Proposition \ref{DegThm}, such cancellations can occur only at the royal nodes $\sigma_j \in \mathbb{T}, j = 1, ..., k$, and then only when $\omega = \overline{x_2(\sigma_j)} = \overline{{\tilde {\eta}}_j}, j = 1, ..., k  $. Hence there are no cancellations in equation (\ref{Cancellation}), and so deg$(\psi_{\omega})=n$. 

(2) Because Problem \ref{The Blaschke interpolation problem} is solvable, its Pick matrix is positive definite and so, by Theorem \ref{prop3.9},  there are polynomials $ a, b, c, d $ of degree at most $n$ that  parametrize the solutions of Problem \ref{The Blaschke interpolation problem}. Choose four particular such  polynomials, as described in Theorem \ref{prop3.9}. By Theorem \ref{prop3.8}, there is $ \tau \in \mathbb{T} \setminus \{ \sigma_{1}, ... , \sigma_{k}\}$ such that the exceptional set $ Z_\tau$ for Problem \ref{The Blaschke interpolation problem}, which is defined as
\begin{equation} \label{ZTauu}
Z_\tau = \{ \zeta \in \mathbb{T} : \text{ for some } j \mbox{ such that } 1\leq j \leq k,  \langle M^{-1}u_{\zeta,\tau} , e_j \rangle =   0 \},
\end{equation} 
consists of at most $k$ points. Fix  a $ \tau \in \mathbb{T} \setminus \{ \sigma_{1}, ... , \sigma_{k}\}$ such that $Z_\tau$ consists of at most $k$ points;  then there exist unique polynomials $ a_{\tau}, b_{\tau}, c_{\tau}, d_{\tau} $ of degree at most $n$ such that 
\begin{equation} \label{equ_3.11_2}
\begin{bmatrix}
a_{\tau}(\tau)       & b_{\tau}(\tau) \\
c_{\tau}(\tau)       & d_{\tau}(\tau) \\
\end{bmatrix}
=
\begin{bmatrix}
1 & 0  \\
0 & 1 \\
\end{bmatrix}
\end{equation}
and, for all $\zeta \in \mathbb{T}  \setminus  Z_\tau $, the function
\begin{equation} \label{ParamEquation}
\varphi = \frac { a_{\tau}\zeta + b_{\tau}} { c_{\tau}\zeta + d_{\tau}}
\end{equation}
is the unique solution of Problem \ref{The Blaschke interpolation problem}  satisfying $ \varphi(\tau) = \zeta$. In addition, the general 4-tuple of polynomials that parametrizes the solutions of Problem \ref{The Blaschke interpolation problem} can be written in the form 
\begin{equation} \label{4Tuple}
(a, b, c, d) = (Xa_{\tau}, Xb_{\tau}, Xc_{\tau}, Xd_{\tau})
\end{equation}
for some rational function $X$.

(3) For $\tau \in \mathbb{T} \setminus  \{ \sigma_{1} , ... , \sigma_{k} \}$ as above, let  $ x_1^{\circ} = x_{1}(\tau) ,  x_2^{\circ} = x_{2}(\tau) ,  x_3^{\circ}= x_{3}(\tau) .$ Since $x$ is tetra-inner, by Theorem \ref{PropOfTetra}, $ |x_3^{\circ}| = 1 $ and $x_1^{\circ} = \overline{x_2^{\circ}} x_3^{\circ}.$ Since $\tau$ is chosen not to be a royal node of $x, |x_1^{\circ}| < 1 ,  |x_2^{\circ}| < 1 $. Thus the relations (\ref{11}) and (\ref{22}) hold. 

\begin{lem} \label{SuppLemma}  
Let $  x_1^{\circ}, x_2^{\circ}, x_3^{\circ} \in \mathbb{C} $ be  such that
\begin{equation} 
| x_3^{\circ} | = 1, \ \ \ | x_1^{\circ} | < 1, \ \ \  | x_2^{\circ} | < 1,
\end{equation}
\begin{equation} 
x_1^{\circ} = \overline{x_2^{\circ}} x_3^{\circ}.
\end{equation} Let $Z_\tau$ be defined as in equation  \eqref{ZTauu}, let $ \tau \in \mathbb{T} \setminus \{ \sigma_{1}, ... , \sigma_{k}\}$ be such that $Z_\tau$ consists of at most $k$ points, and let 					
	$$
	{Z ^{\sim}_\tau} = \Big\{ {\frac {{\overline{{\tilde {\eta}}_1}} x_3^{\circ} - x_1^{\circ}} {x_2^{\circ} {\overline{{\tilde {\eta}}_1}} - 1 }} , {\frac { {\overline{{\tilde {\eta}}_2}} x_3^{\circ} - x_1^{\circ} } {x_2^{\circ} {\overline{{\tilde {\eta}}_2}} - 1 }   } , . . . , {\frac { {\overline{{\tilde {\eta}}_k}} x_3^{\circ} - x_1^{\circ} } {x_2^{\circ} {\overline{{\tilde {\eta}}_k}} - 1 }   } \Big\}.              
	$$						
	If $\zeta \in \mathbb{T}  \setminus  {Z^{\sim}_\tau} $ then the function
	\begin{equation} \label{Varsolution}
	\varphi = \frac { (x_2^{\circ} x_1 - x_3 )\zeta + x_1^{\circ} x_3 - x_1 x_3^{\circ}} { (x_2^{\circ} - x_2)\zeta + x_1^{\circ} x_2 - x_3^{\circ} }
	\end{equation}		
	is a solution of {\em Problem \ref{The Blaschke interpolation problem} } which satisfies $ \varphi(\tau) = \zeta $. 
\end{lem}
\begin{proof} By equation \eqref{key equation}, for any $ \omega \in \mathbb{T}$,	
	$$\psi_{\omega}(\tau) = \frac { \omega x_3^{\circ} - x_1^{\circ} }  { x_2^{\circ}  \omega - 1  } ,$$
which is well defined since $| x_{2}^{\circ} | < 1$. We have, for $\zeta \in \mathbb{T}, $
	$$
	\psi_{\omega}(\tau) = 
\zeta \Leftrightarrow \frac { \omega x_3^{\circ} - x_1^{\circ} }  { x_2^{\circ}  \omega - 1  }  = \zeta \Leftrightarrow \omega = \frac {-\zeta + x_1^{\circ} } {x_3^{\circ} - \zeta x_2^{\circ}}.
	$$
	Hence, as long as 
	\begin{equation} \label{nonequal}
	\frac {-\zeta + x_1^{\circ} } {x_3^{\circ} - \zeta x_2^{\circ}} \neq  \ {\overline{{\tilde {\eta}}_1}} , .... , {\overline{{\tilde {\eta}}_k}},
	\end{equation}
	the function
	\begin{eqnarray*}
		\varphi(\lambda) = \psi_{\omega}(\lambda) &=& \psi_{\frac {-\zeta + x_1^{\circ} } {x_3^{\circ} - \zeta x_2^{\circ}}} (\lambda) \\
		&=&  \frac{   \displaystyle{\frac{x_{3}(\lambda)x_1^{\circ} - x_{3}(\lambda) \zeta }{x_3^{\circ} - x_2^{\circ} \zeta} } -  x_{1}(\lambda) }{ \displaystyle{\frac{x_{2}(\lambda) x_1^{\circ} - x_{2}(\lambda) \zeta }{x_3^{\circ} - x_2^{\circ} \zeta} } - 1         } \\
&=& \frac { ( x_1(\lambda) x_2^{\circ} - x_3(\lambda) )\zeta + x_1^{\circ} x_3(\lambda) - x_1(\lambda) x_3^{\circ}} { (x_2^{\circ} - x_2(\lambda))\zeta + x_1^{\circ} x_2(\lambda) - x_3^{\circ} }.
	\end{eqnarray*}
	is a solution of Problem \ref{The Blaschke interpolation problem} which satisfies  $ \varphi(\tau) = \zeta$.
	Condition (\ref{nonequal}) can also be written, for $ j = 1,2,...,k, $ 
	$$
	\zeta \neq \frac{x_1^{\circ} - {\overline{{\tilde {\eta}}_j}} x_3 ^{\circ}  } { 1 - x_2^{\circ} {\overline{{\tilde {\eta}}_j}}   } = \frac{{\overline{{\tilde {\eta}}_j}} x_3 ^{\circ} - x_1^{\circ} }{x_2^{\circ} {\overline{{\tilde {\eta}}_j}} - 1}
	$$
	or equivalently $\zeta \notin {Z^{\sim}_\tau}. $ 
\end{proof}
For $\zeta \in \mathbb{T}  \setminus  ( Z_{\tau} \cup {Z^{\sim}_\tau} )$ where ${Z^{\sim}_\tau}$ is defined in Lemma \ref{SuppLemma}, we have two formulae for the unique solution of Problem \ref{The Blaschke interpolation problem} satisfying $ \varphi(\tau) = \zeta $, namely, the equations (\ref{ParamEquation})  and (\ref{Varsolution}). Note that 
\begin{eqnarray*}
\begin{bmatrix}
x_2^{\circ} x_{1}(\tau) - x_{3}(\tau)       &  x_1^{\circ} x_{3}(\tau) - x_{1}(\tau) x_3^{\circ}  \\
x_2^{\circ} - x_{2}(\tau)       & x_1^{\circ} x_{2}(\tau) - x_3^{\circ}  \\
\end{bmatrix}
&=&
\begin{bmatrix}
x_2^{\circ} x_1^{\circ} - x_3^{\circ}       &  x_1^{\circ} x_3^{\circ} - x_1^{\circ} x_3^{\circ}  \\
x_2^{\circ} - x_2^{\circ}       & x_1^{\circ} x_2^{\circ} - x_3^{\circ}  \\
\end{bmatrix} \\ 
&=&
(x_1^{\circ} x_2^{\circ} - x_3^{\circ}) \begin{bmatrix}
1 & 0  \\
0 & 1 \\
\end{bmatrix}.
\end{eqnarray*}
Because the set $Z_{\tau} \cup {Z^{\sim}_\tau}$ is finite, the linear fractional transformations in equations (\ref{ParamEquation}) and (\ref{Varsolution}) are equal at infinitely many points and therefore coincide. It follows from the normalizing condition that
\begin{equation}
\begin{bmatrix}
a_{\tau}     & b_{\tau} \\
c_{\tau}      & d_{\tau} \\
\end{bmatrix}
=
\frac{1}{x_1^{\circ} x_2^{\circ} - x_3^{\circ}}
\begin{bmatrix}
x_2^{\circ} x_{1} - x_{3}       &  x_1^{\circ} x_{3} - x_{1} x_3^{\circ}  \\
x_2^{\circ} - x_{2}       & x_1^{\circ} x_{2} - x_3^{\circ}  \\
\end{bmatrix}.
\end{equation}
Suppose that $a, b, c \ \text{and} \  d $ are polynomials that parametrize the solutions of Problem \ref{The Blaschke interpolation problem}, as in Theorem \ref{4.4Tetra} (2). By the observation (\ref{4Tuple}), there is a rational function $X$ such that
\begin{equation}\label{Xa}
Xa = x_2^{\circ} x_1 - x_3,
\end{equation}
\begin{equation}\label{Xb}
Xb = x_1^{\circ} x_3 - x_1 x_3^{\circ},
\end{equation}
\begin{equation} \label{Xc}
Xc = x_2^{\circ} - x_2,
\end{equation}
\begin{equation} \label{Xd}
Xd = x_1^{\circ} x_2 - x_3^{\circ},
\end{equation}
Let us find connections between $x_1, x_2,x_3$ and the polynomials $ a, b, c, d $. Equations (\ref{Xc}) and (\ref{Xd}) for $x_2$ and $X$ can  be written as
\begin{equation} \label{linearsys1}
\begin{array}{ccl}
Xc + x_2 & = & x_2^{\circ} \\
 Xd - x_1^{\circ} x_2   & = & - x_3^{\circ}. 
\end{array}
\end{equation}
Then, the solution of the system (\ref{linearsys1}) is
\begin{equation} \label{X}
X = \frac{\begin{vmatrix} x_2^{\circ} & 1 \\ -x_3^{\circ} & -x_1^{\circ} \end{vmatrix}}{\begin{vmatrix} c & 1 \\ d & -x_1^{\circ} \end{vmatrix}} = \frac{x_1^{\circ} x_2^{\circ} - x_3^{\circ}}{x_1^{\circ} c +d}
\end{equation}
and
\begin{equation}
x_2 = \frac{\begin{vmatrix} c & x_2^{\circ} \\ d & -x_3^{\circ} \end{vmatrix}}{\begin{vmatrix} c & 1 \\ d & -x_1^{\circ} \end{vmatrix}} = \frac{x_3^{\circ} c + x_2^{\circ} d } { x_1^{\circ} c + d }.
\end{equation}
Equations (\ref{Xa}) and (\ref{Xb}) give us the system
\begin{equation} \label{linearsys2}
\begin{array}{ccl}
  x_2^{\circ} x_1 - x_3 & = & Xa \\
	-x_3^{\circ} x_1 +  x_1^{\circ} x_3   & = & Xb . 
\end{array}
\end{equation}
Then, the solution of the system (\ref{linearsys2}) is 
\begin{eqnarray*}
x_1 = \frac{\begin{vmatrix} Xa & -1 \\ Xb & x_1^{\circ} \end{vmatrix}}{\begin{vmatrix} x_2^{\circ} & -1 \\ -x_3^{\circ} & x_1^{\circ} \end{vmatrix}}
&=& \frac{\begin{vmatrix} \displaystyle{\frac{x_1^{\circ} x_2^{\circ}a - x_3^{\circ}a}{x_1^{\circ} c +d}} & -1 \\ \displaystyle{\frac{x_1^{\circ} x_2^{\circ}b - x_3^{\circ}b}{x_1^{\circ} c +d}} & x_1^{\circ} \end{vmatrix}}{\begin{vmatrix} x_2^{\circ} & -1 \\ -x_3^{\circ} & x_1^{\circ} \end{vmatrix}} \\
&=& \frac { x_1^{\circ} a + b } { x_1^{\circ} c + d }
\end{eqnarray*} 
and
\begin{eqnarray*}
	x_3 = \frac{\begin{vmatrix}  x_2^{\circ} & Xa \\ -x_3^{\circ} & Xb \end{vmatrix}}{\begin{vmatrix} x_2^{\circ} & -1 \\ -x_3^{\circ} & x_1^{\circ} \end{vmatrix}}	
	&=& \frac{\begin{vmatrix}  x_2^{\circ} & \displaystyle{\frac{x_1^{\circ} x_2^{\circ}a - x_3^{\circ}a}{x_1^{\circ} c +d}}  \\ -x_3^{\circ} & \displaystyle{\frac{x_1^{\circ} x_2^{\circ}b - x_3^{\circ}b}{x_1^{\circ} c +d}} \end{vmatrix}}{\begin{vmatrix} x_2^{\circ} & -1 \\ -x_3^{\circ} & x_1^{\circ} \end{vmatrix}} \\ 
	&=& \frac { x_2^{\circ} b + x_3^{\circ} a } { x_1^{\circ} c + d }.
\end{eqnarray*}
Thus $x_1, x_2, x_3$ are given by equations (\ref{x_1Equation}), (\ref{x_2Equation}) and (\ref{x_3Equation}). The proof of Theorem \ref{4.4Tetra} is complete. 
\end{proof}

Note that we can also prove a result similar to Theorem \ref{4.4Tetra}, using the function $\Upsilon_\omega$ instead of $\Psi_{\omega}$, where 
$$\displaystyle{ \Upsilon_\omega(x_1, x_2, x_3) =  \frac{x_{3} \omega -x_{2}}{x_{1}\omega-1}}, $$ which is defined for every $( x_1, x_2, x_3)$ in $\mathbb{C}^{3}$ such that $ x_{1} \omega - 1 \neq 0 $. In this case we suppose that $\rho_j =  A x_2  (\sigma_j)$ for $ j = 1, 2 , .. , k .$

\section{From the Blaschke interpolation problem to the royal tetra-interpolation problem} \label{FromBlashtoRoyal}
In this section we will prove Theorem \ref{ReverseThm}. This theorem shows that, if Blaschke interpolation  data $(\sigma, \eta, \rho)$  of type $(n,k)$ are given and the corresponding Problem \ref{The Blaschke interpolation problem} is solvable, then we can construct a solution for the royal tetra-interpolation problem $(\sigma, \eta, \tilde{\eta}, \rho)$, for some $\tilde{\eta} = (\tilde{\eta_1}, ..., \tilde{\eta_n})$. We will start with some technical lemmas. 
\begin{lem} \label{x_2Lemma}
	Let $a, b, c, d, x_1^{\circ}, x_2^{\circ}, x_3^{\circ} \in \mathbb{C}$, and suppose that $|x_3^{\circ}| = 1, x_1^{\circ} = \overline{x_2^{\circ}} x_3^{\circ}, x_1^{\circ} c \neq -d $ and $|x_1^{\circ}| < 1, |x_2^{\circ}| < 1$. Let 
	$$
	x_2 = \frac{x_3^{\circ} c + x_2^{\circ} d}{x_1^{\circ} c + d}.
	$$
	Then
	$$
	(1) \ \ |x_2| \leq 1 \ \text{ if and only if} \ |c| \leq |d|,
	$$
	and
	$$
	(2) \ \ |x_2| < 1 \ \text{ if and only if} \ |c| < |d|.
	$$
\end{lem}
	\begin{proof} (1)
		\begin{eqnarray}\label{cleqd} \nonumber
	|x_2| \leq 1 \ 	&\Leftrightarrow& \  {\Bigg|\frac{x_3^{\circ} c + x_2^{\circ} d}{x_1^{\circ} c + d} \Bigg|} \leq 1 \\  \nonumber
	&\Leftrightarrow& |x_3^{\circ} c + x_2^{\circ} d|^{2} \leq |x_1^{\circ} c + d|^{2} \\  \nonumber
	&\Leftrightarrow& |c|^{2} + 2\text{Re}(x_3^{\circ} \ c \ \overline{x_2^{\circ}} \ \overline{d}) +| x_2^{\circ}|^{2} |d|^{2} \leq |x_1^{\circ}|^{2} |c|^{2} + 2\text{Re} (x_1^{\circ} \ c \ \overline{d}) + |d|^{2}  \\  \nonumber 
	&\Leftrightarrow& |c|^{2} + |x_2|^{2} |d|^{2} - |x_1^{\circ}|^{2} |c|^{2} - |d|^{2} \leq 0 \ \ (\text{since} \ x_1^{\circ} = \overline{x_2^{\circ}} x_3^{\circ})\\  \nonumber
	&\Leftrightarrow& (1-|x_1^{\circ}|^{2}) (|c|^{2} - |d|^{2}) \leq 0 \\  
	&\Leftrightarrow& |c| \leq |d| \quad (\text{since} \  (1-|x_1^{\circ}|^{2}) > 0). 
		\end{eqnarray}
		
		(2) The same calculation leads to 	$|x_2| < 1 \ \iff\ |c| < |d| $.
	\end{proof}

\begin{prop} \label{AlgeberaicProp}
	Let a, b, c, d be polynomials in the variable $\lambda$ and let $x_1^{\circ}, x_2^{\circ}, x_3^{\circ} \in \mathbb{C}$ satisfy $x_3^{\circ} \neq x_1^{\circ} x_2^{\circ} $ and $ x_1^{\circ} c \neq -d$.
	Let rational functions $x_1, x_2, x_3$ be defined by 
\begin{equation} \label{X1X2X3}
x_1(\lambda) = \frac { x_1^{\circ} a(\lambda) + b(\lambda) } { x_1^{\circ} c(\lambda) + d(\lambda) }, \ 
x_2(\lambda) = \frac{x_3^{\circ} c(\lambda) + x_2^{\circ} d(\lambda) } { x_1^{\circ} c(\lambda) + d(\lambda) }, \ 
x_3(\lambda) = \frac { x_2^{\circ} b(\lambda) + x_3^{\circ} a(\lambda)} { x_1^{\circ} c(\lambda) + d(\lambda) }.
\end{equation}
and define a rational function $\zeta$ in the indeterminate $\omega$ by
\begin{equation}
\zeta(\omega) = \frac { \omega x_3^{\circ} - x_1^{\circ} }  { x_2^{\circ}  \omega - 1  } .
\end{equation}
Then, as rational functions in $(\omega, \lambda)$,
$$
\frac { \omega x_3(\lambda) - x_1(\lambda) }  { x_2(\lambda)  \omega - 1  } = \frac{a(\lambda)\zeta(\omega) + b(\lambda)}{c(\lambda) \zeta(\omega) + d(\lambda)}.
$$
\end{prop}


\begin{proof}
Let  $x_1 , x_2, x_3$ be defined by equations (\ref{X1X2X3}).
Then
\begin{eqnarray*}
\displaystyle{\frac { \omega x_3(\lambda) - x_1(\lambda) }  { x_2(\lambda)  \omega - 1  }} &=&\frac{ \displaystyle{\frac { \omega x_2^{\circ} b(\lambda) + \omega x_3^{\circ} a(\lambda)} { x_1^{\circ} c(\lambda) + d(\lambda)} - \displaystyle{\frac { x_1^{\circ} a(\lambda) + b(\lambda) } { x_1^{\circ} c(\lambda) + d(\lambda) }} }}{ \displaystyle{\frac{ \omega x_3^{\circ} c(\lambda) + \omega x_2^{\circ} d(\lambda) } { x_1^{\circ} c(\lambda) + d(\lambda) }} -1} \\ &=& \frac{\omega x_2^{\circ} b(\lambda) + \omega x_3^{\circ} a(\lambda) - x_1^{\circ} a(\lambda) - b(\lambda)}{\omega x_3^{\circ} c(\lambda) + \omega x_2^{\circ} d(\lambda) - x_1^{\circ} c(\lambda) - d(\lambda)} \\
&=& \frac{a(\lambda) (\omega x_3^{\circ} - x_1^{\circ} ) + b(\lambda) (\omega x_2^{\circ} - 1)}{c(\lambda) (\omega x_3^{\circ} - x_1^{\circ}) + d(\lambda) (\omega x_2^{\circ} - 1)} \\
&=& \frac{a(\lambda) \left(\displaystyle{\frac{\omega x_3^{\circ} - x_1^{\circ} }{\omega x_2^{\circ} - 1}}\right) + b(\lambda)}{c(\lambda) \left(\displaystyle{\frac{\omega x_3^{\circ} - x_1^{\circ} }{\omega x_2^{\circ} - 1}} \right) + d(\lambda )} \\
&=& \frac{a(\lambda) \zeta(\omega) + b(\lambda)}{c(\lambda) \zeta(\omega) + d(\lambda)},
\end{eqnarray*}
where $\zeta(\omega) = \displaystyle{\frac{\omega x_3^{\circ} - x_1^{\circ} }{\omega x_2^{\circ} - 1} }$.
\end{proof}
\begin{prop} \label{Prooop} 
	Let a, b, c, d be polynomials having no common zero in $\overline{\mathbb{D}}$, and satisfying $|c| \leq |d|$ on $\mathbb{D}$ . Suppose that $x_1^{\circ}, x_2^{\circ}, x_3^{\circ} \in \mathbb{C}$ satisfy $x_1^{\circ} c \neq -d$,  $|x_3^{\circ}| = 1, |x_1^{\circ}| < 1 , |x_2^{\circ}| < 1$ and $x_1^{\circ} = \overline{x_2^{\circ}} x_3^{\circ}$. Let rational functions $x_1, x_2, x_3$ be defined by
	\begin{equation} \label{x'sequations}
	x_1(\lambda) = \frac { x_1^{\circ} a(\lambda) + b(\lambda) } { x_1^{\circ} c(\lambda) + d(\lambda) }, \ 
	x_2(\lambda) = \frac{x_3^{\circ} c(\lambda) + x_2^{\circ} d(\lambda) } { x_1^{\circ} c(\lambda) + d(\lambda) }, \ 
	x_3(\lambda) = \frac { x_2^{\circ} b(\lambda) + x_3^{\circ} a(\lambda)} { x_1^{\circ} c(\lambda) + d(\lambda) },
	\end{equation}
	and let 
	\begin{equation} \label{uppsiRelation}
	\psi_\zeta(\lambda) = \frac{a(\lambda)\zeta + b(\lambda)}{c(\lambda) \zeta + d(\lambda)}.
	\end{equation}
	\begin{enumerate}
		\item[(i)] If, for all but finitely many values of $\lambda \in \mathbb{D}$,
		\begin{equation} \label{inD}
		|	\psi_\zeta(\lambda)| \leq 1
		\end{equation}
for all but finitely many $\zeta \in \mathbb{T}$, then $x_1^{\circ} c +d$ has no zero in $\overline{\mathbb{D}}$ and $ x = (x_1, x_2, x_3)$ is an analytic map from $\mathbb{D} $ to $\overline{\mathcal{E}}$.
		\item[(ii)] If, for all but finitely many $\zeta \in \mathbb{T}$, the function $\psi_\zeta$ is inner, then either $x(\overline{\mathbb{D}}) \subseteq \mathcal{R}_{\overline{\mathcal{E}}} $ or $ x = (x_1, x_2, x_3)$ is a rational tetra-inner function.
 	\end{enumerate}
\end{prop}
\begin{proof}
	(i) Suppose there is a finite subset $E$ of $\mathbb{D}$ such that, for all $\lambda \in \mathbb{D} \setminus E$, there is a finite subset $F_\lambda$ of $\mathbb{T}$ for which the inequality (\ref{inD}) holds for all $\zeta \in \mathbb{T} \setminus F_\lambda$. 
	Let us show that the denominator  $x_1^{\circ} c +d$ of $x_1, x_2, x_3$ has no zeros in $\overline{\mathbb{D}}$. Suppose that $\alpha \in \overline{\mathbb{D}}$ is a zero of $(x_1^{\circ}  c + d)$. Since $|c| \leq |d|$ on $\overline{\mathbb{D}}$,
	\begin{eqnarray*}
 |x_1^{\circ} c(\alpha) + d(\alpha)|  &\geq& |d(\alpha)| - |x_1^{\circ}  c(\alpha) |\\
	&\geq& |d(\alpha)| - |x_1^{\circ} | |d(\alpha)| \\
	&=& (1-|x_1^{\circ} |) |d(\alpha)| .
\end{eqnarray*}
Thus,
$$
0 = |x_1^{\circ} c(\alpha) + d(\alpha)| \geq (1-|x_1^{\circ} |) |d(\alpha)|.
$$
Since $|x_1^{\circ}| < 1, \ (1-|x_1^{\circ} |) \neq 0 $, and so $ d(\alpha) = 0 $, Then
$$
0 = x_1^{\circ} c(\alpha) + d(\alpha) = x_1^{\circ} c(\alpha)
$$
implies that  $c(\alpha) = 0 $.
	
	Pick a sequence $(\alpha_j)$ in $\mathbb{D} \setminus E$ such that $\alpha_j \rightarrow \alpha$. For each $j$, for $\zeta \in \mathbb{T} \setminus F_{\lambda_j}$, we have $|	\psi_\zeta(\lambda)| \leq 1$ on $\mathbb{D} \setminus E$. Hence, for $\zeta \in \mathbb{T} \setminus \cup_j F_{\lambda_j} $, which is to say, for all but countably many $\zeta\in\T$,
	$$
	\bigg|\frac{a(\alpha_j)\zeta + b(\alpha_j)}{c(\alpha_j) \zeta + d(\alpha_j)} \bigg| \leq 1. 
	$$
Because $c(\alpha_j) \zeta + d(\alpha_j) \rightarrow 0$ uniformly almost everywhere for $\zeta \in \mathbb{T} $ as $j \rightarrow \infty$, the same holds for $a(\alpha_j)\zeta + b(\alpha_j)$. Therefore $ a(\alpha_j) \rightarrow 0 $ and $b(\alpha_j) \rightarrow 0 $. Hence $a(\alpha) = b(\alpha) = 0$.
Thus $a, b, c, d$ all vanish at $\alpha$, contrary to our assumption. So $x_1^{\circ} c +d$ has no zeros in $\overline{\mathbb{D}}$. Thus $x_1, x_2, x_3$  defined by equations (\ref{x'sequations}) are rational functions having no poles in $\overline{\mathbb{D}}$.
	
	Consider $\lambda \in \mathbb{D} \setminus E$. By Proposition \ref{AlgeberaicProp},
	\begin{equation}\label{UppsiRelation}
	\Psi_{\omega}(x_1(\lambda), x_2(\lambda), x_3(\lambda)) = \frac { \omega x_3(\lambda) - x_1(\lambda) }  { x_2(\lambda)  \omega - 1  } = \frac{a(\lambda)\zeta(\omega) + b(\lambda)}{c(\lambda) \zeta({\omega}) + d(\lambda)} 
	\end{equation}
provided that both sides are defined, that is, for all $\omega \in \mathbb{T} \setminus \Omega_\lambda $ where
$$
\Omega_\lambda = \{ \omega \in \mathbb{T} : \omega x_2({\lambda}) =1 \ \text{or} \ c(\lambda) \zeta(\omega) = -d(\lambda)  \}.
$$	
$\Omega_\lambda$ contains at most two points. On combining the relations (\ref{uppsiRelation}), (\ref{inD}) and (\ref{UppsiRelation}), we deduce that, for $\lambda \in \mathbb{D} \setminus E,$
\begin{equation} \label{WW}
|\Psi_{\omega}(x_1(\lambda), x_2(\lambda), x_3(\lambda))| \leq 1
\end{equation}
for all $\omega \in \mathbb{T}$ such that $\omega \notin \Omega_\lambda \cup \zeta^{-1} (F_\lambda)$, that is, for all but finitely many $\omega \in \mathbb{T}$. By \cite[Theorem 2.4]{AWY} (see Theorem \ref{DefOfTetra}), $(x_1(\lambda), x_2(\lambda), x_3(\lambda)) \in \overline{\mathcal{E}}$. Since this is true for all but finitely many   $\lambda \in \mathbb{D}$, and $x_1, x_2, x_3$ are rational functions without poles in $\overline{ \mathbb{D}}$, $(x_1, x_2, x_3)$ maps $\mathbb{D}$ into $\overline{\mathcal{E}}$.

(ii) Let us assume that, for some finite subset $F$ of $\mathbb{T}$, the function $\psi_\zeta$ is inner for all $\zeta \in \mathbb{T} \setminus F$. By Part (i), $(x_1, x_2, x_3)$ is a rational analytic map from $\mathbb{D}$ into $\overline{\mathcal{{E}}}$ and therefore extends to a continuous map of $\overline{\mathbb{D}}$ into $\overline{\mathcal{{E}}}$. Let $\lambda \in \mathbb{T}$. By Proposition \ref{AlgeberaicProp} and equation (\ref{uppsiRelation}),
\begin{equation}
\Psi_{\omega}(x_1(\lambda), x_2(\lambda), x_3(\lambda)) = \psi_{\zeta(\omega)}(\lambda) 
\end{equation}
provided that both sides are defined, that is, for all $\omega \in \mathbb{T} \setminus \Omega_\lambda $ where
$$
\Omega_\lambda = \{ \omega \in \mathbb{T} : \omega x_2({\lambda}) =1 \ \text{or} \ c(\lambda) \zeta(\omega) = -d(\lambda)  \}.
$$	
Note that $\Omega_\lambda$ contains at most two points. For $\omega \in \mathbb{T} \setminus \zeta^{-1} (F)$ the function $\psi_{\zeta(\omega)}$ is inner.
Thus, for $\omega \in \mathbb{T} \setminus (\zeta^{-1} (F) \cup \Omega_\lambda)$,
\begin{equation} \label{equalityEquation}
|\Psi_{\omega}(x_1(\lambda), x_2(\lambda), x_3(\lambda))| = |\psi_{\zeta(\omega)}(\lambda) | = 1 .
\end{equation}
Case 1. Suppose that for all $\lambda \in \overline{ \mathbb{D}}$,  $x_1(\lambda) x_2(\lambda) = x_3(\lambda)$. Then, for all $\lambda\in \overline{\mathbb{D}}$,
\begin{align*}
\Psi_{\omega}(x_1(\lambda), x_2(\lambda), x_3(\lambda)) = \frac { \omega x_3(\lambda) - x_1(\lambda) }  { x_2(\lambda)  \omega - 1  }  &=\frac{\omega x_1(\lambda) x_2(\lambda) - x_1(\lambda) }{x_2(\lambda)  \omega - 1} \\
&=   \frac{x_1(\lambda) (\omega x_2(\lambda) - 1)}{x_2(\lambda)  \omega - 1} = x_1(\lambda) .  \\
\end{align*} Thus  $x(\overline{\mathbb{D}}) \subseteq \mathcal{R}_{\overline{\mathcal{E}}} $. 
\newline Case 2.  Suppose that for some $\lambda \in \overline{ \mathbb{D}}$,  $x_1(\lambda) x_2(\lambda) \neq x_3(\lambda)$.  To prove that $x=(x_1, x_2,x_3)$ is rational $\overline{\mathcal{{E}}}$-inner function, by Theorem \ref{PropOfTetra}, we need to show that   
 $(x_1, x_2,x_3)(\lambda) \in b\overline{\mathcal{{E}}}$ for almost all $\lambda \in \mathbb{T}$, that is, 
\begin{enumerate}
 	\item $ |x_3(\lambda)| =1$ for almost all $\lambda \in \mathbb{T}$,
 	\item $|x_2| \leq 1 $ on $\overline{\mathbb{D}}$, 
 	\item $x_1(\lambda) = \overline{x_2(\lambda)} x_3(\lambda) $ for almost all $ \lambda \in \mathbb{T}$.
 \end{enumerate}
For (2), by Lemma \ref{x_2Lemma}, we showed that  $|x_2(\lambda)| \leq 1 $ for $\lambda \in \mathbb{\overline{D}}$. By Lemma \ref{InnerUppsi}, for any  $\omega \in \mathbb{T}$ and any point $x = (x_1, x_2, x_3) \in \mathcal{\overline{E}}$ such that $x_1 x_2 \neq x_3$, 
$$
|\Psi_{\omega}(x_1, x_2, x_3)| = 1 \ \text{if and only if } \  2 \omega (x_2 - \overline{x_1} x_3 )) = 1  - |x_1|^{2} + |x_2|^{2} - |x_3|^{2}.
$$
Thus, for $\lambda \in \mathbb{T}$ such that $x_1(\lambda) x_2(\lambda) \neq x_3(\lambda)$ , equation (\ref{equalityEquation}) implies
$$
2 \omega (x_2(\lambda) - \overline{x_1(\lambda)} x_3(\lambda) )) = 1  - |x_1(\lambda)|^{2} + |x_2(\lambda)|^{2} - |x_3(\lambda)|^{2}.
$$
Hence, for $\lambda \in \mathbb{T}$, if $|\Psi_{\omega}(x_1(\lambda), x_2(\lambda), x_3(\lambda))| = 1$ for two distinct $\omega \in \mathbb{T}$, 
say
 $\omega_1 \neq \omega_2$, we have the linear system
\begin{equation} \label{Equality.}
\begin{array}{ccl}
2\omega_1 (x_2(\lambda) - \overline{x_1(\lambda)} x_3(\lambda) ) & = & 1  - |x_1(\lambda)|^{2} + |x_2(\lambda)|^{2} - |x_3(\lambda)|^{2}  \\
2 \omega_2 (x_2(\lambda) - \overline{x_1(\lambda)} x_3(\lambda) )   & = & 1  - |x_1(\lambda)|^{2} + |x_2(\lambda)|^{2} - |x_3(\lambda)|^{2}.
\end{array}
\end{equation}
Thus, for $\lambda \in \mathbb{T}$,
$$
2 \omega_1 (x_2(\lambda) - \overline{x_1(\lambda)} x_3(\lambda) ) - 2 \omega_2 (x_2(\lambda) - \overline{x_1(\lambda)} x_3(\lambda) ) = 0
$$
\begin{align} \nonumber
&\implies (x_2(\lambda) - \overline{x_1(\lambda)} x_3(\lambda) ) ( \omega_1 - \omega_2 ) = 0\\ 
&\implies x_2(\lambda) = \overline{x_1(\lambda)} x_3(\lambda). 
\end{align}
By equations (\ref{Equality.}), for $\lambda \in \mathbb{T}$,  
\begin{equation}\label{equalitytoZero}
1  - |x_1(\lambda)|^{2} + |x_2(\lambda)|^{2} - |x_3(\lambda)|^{2} = 0 .
\end{equation}
Note for $\lambda \in \mathbb{T}$, since $x_2(\lambda) = \overline{x_1(\lambda)} x_3(\lambda)$, 
\begin{align} \nonumber
(\ref{equalitytoZero}) \ \mbox{holds} \
&\Leftrightarrow 1  - |x_1(\lambda)|^{2} + |\overline{x_1(\lambda)} x_3(\lambda)|^{2} - |x_3(\lambda)|^{2} = 0 \\ \nonumber
&\Leftrightarrow 1  - |x_1(\lambda)|^{2} + |\overline{x_1(\lambda)}|^{2} | x_3(\lambda)|^{2} - |x_3(\lambda)|^{2} = 0 \\ \nonumber
&\Leftrightarrow 1  - |x_1(\lambda)|^{2} - | x_3(\lambda)|^{2} (1 - |{x_1(\lambda)}|^{2} ) = 0 \\ \nonumber
&\Leftrightarrow   (1  - |x_1(\lambda)|^{2}) (1- |x_3(\lambda)|^{2}) = 0 \\ \nonumber
&\Leftrightarrow |x_3(\lambda)| = 1 \ \text{or} \ |x_1(\lambda)| = 1.
\end{align} \ 
Case $1$. If $|x_1(\lambda)| = 1$ and $x_2(\lambda) = \overline{x_1(\lambda)} x_3(\lambda)$, we have $x_3(\lambda)= x_1(\lambda) x_2(\lambda)  $ for almost all $\lambda \in \mathbb{T}$. Then since $x_i$ are rational functions for $i = 1,2,3$, and $x_3(\lambda)= x_1(\lambda) x_2(\lambda)  $ for $\lambda \in \mathbb{T}$, it imples that
$$
 x_3(\lambda) = x_1(\lambda) x_2(\lambda) \quad \text{for all} \ \lambda\in \overline{\mathbb{D}}.
$$ 
Then, for all $\lambda\in \overline{\mathbb{D}}$,
\begin{align*}
\Psi_{\omega}(x_1(\lambda), x_2(\lambda), x_3(\lambda)) = \frac { \omega x_3(\lambda) - x_1(\lambda) }  { x_2(\lambda)  \omega - 1  }  &=\frac{\omega x_1(\lambda) x_2(\lambda) - x_1(\lambda) }{x_2(\lambda)  \omega - 1} \\
&=   \frac{x_1(\lambda) (\omega x_2(\lambda) - 1)}{x_2(\lambda)  \omega - 1} = x_1(\lambda) .  \\
\end{align*} Thus  $x(\overline{\mathbb{D}}) \subseteq \mathcal{R}_{\overline{\mathcal{E}}} $. 
 
Case $2$. If for almost all  $\lambda \in \mathbb{T} , | x_3(\lambda)| = 1$, then 
\begin{align*}
x_2(\lambda) = \overline{x_1(\lambda)} x_3(\lambda) \ &\implies \ x_2(\lambda) \overline{x_3(\lambda)} = \overline{x_1(\lambda)} \\
&\implies \  x_1(\lambda) = \overline{x_2(\lambda)} x_3(\lambda).  \\
\end{align*}
Thus, for almost all $\lambda \in \mathbb{T}$,  $| x_3(\lambda)| = 1$ and $x_1(\lambda) = \overline{x_2(\lambda)} x_3(\lambda)$ that proves (1) and (3) respectively. Therefore,  the point $(x_1(\lambda), x_2(\lambda), x_3(\lambda)) $ for almost all $\lambda \in \mathbb{T}
$ is in the distinguished boundary $b\mathcal{\overline{E}}$ of $\mathcal{\overline{E}}$. Hence $x = (x_1, x_2, x_3)$ is a rational $\mathcal{\overline{E}}$-inner function in this case.  
\end{proof}

\begin{thm} \label{ReverseThm}
	Let $(\sigma, \eta,  \rho)$ be Blaschke interpolation data  of type $(n,k)$, and let $(\sigma, \eta, \tilde{\eta}, \rho)$ be royal tetra-interpolation data of type $(n,k)$ where $\tilde{\eta}= ( \tilde{\eta}_1, \tilde{\eta}_2, ..., \tilde{\eta}_n), $ $\tilde{\eta}_j \in \mathbb{T}, j = 1, ..., k$ and $\tilde{\eta}_j \in \mathbb{D}, j = k+1, ..., n$. Suppose that Problem {\em \ref{The Blaschke interpolation problem}} with $(\sigma, \eta, \rho)$  is solvable and the solutions $\varphi$ of Problem {\em \ref{The Blaschke interpolation problem}} have normalized parametrization 
	$$
	\varphi = \frac { a\zeta + b} { c\zeta + d}.
	$$
	Suppose that there exist scalars $x_1^{\circ}, x_2^{\circ}, x_3^{\circ}$ in $\mathbb{C}$ such that 
	$$
	| x_3^{\circ} | = 1, \ \ \ | x_1^{\circ} | < 1, \ \ \  | x_2^{\circ} | < 1, \ \	x_1^{\circ} = \overline{x_2^{\circ}} x_3^{\circ},
	$$
	and 
	\begin{equation} \label{ExtraCondition}
	\frac{x_3^{\circ} c(\sigma_j) + x_2^{\circ} d(\sigma_j)}{x_1^{\circ} c(\sigma_j) + d(\sigma_j)} = \tilde{\eta_j} , \ j = 1, ..., n. 
	\end{equation}  
	Then there exists a rational tetra-inner function $x = (x_1, x_2, x_3)$  given by
	\begin{equation} \label{Equationforx_1}
	x_1(\lambda) = \frac { x_1^{\circ} a(\lambda) + b(\lambda) } { x_1^{\circ} c(\lambda) + d(\lambda) }
	\end{equation}
	\begin{equation} \label{Equationforx_2}
	x_2(\lambda) = \frac{x_3^{\circ} c(\lambda) + x_2^{\circ} d(\lambda) } { x_1^{\circ} c(\lambda) + d(\lambda) }
	\end{equation}
	\begin{equation} \label{Equationforx_3}
	x_3(\lambda) = \frac { x_2^{\circ} b(\lambda) + x_3^{\circ} a(\lambda)} { x_1^{\circ} c(\lambda) + d(\lambda) },
	\end{equation} for $\lambda \in \mathbb{D}$, such that
	\begin{enumerate}
		\item[(i)] $ x \in \mathcal{R}^{n,k}$, and $x$ is a solution of the royal tetra-interpolation problem with the data $(\sigma, \eta, \tilde{\eta}, \rho)$, that is, 
		$$x(\sigma_j) = (\eta_j, \tilde{\eta_j}, \eta_i \tilde{\eta_j}) \ \text{for} \ j = 1, ... ,n,$$
		and 
		$$Ax_1(\sigma_j) = \rho_j \ \text{for} \  j = 1, ... , k,$$
		\item[(ii)] for all but finitely many $\omega \in \mathbb{T}$, the function $\Psi_{\omega} \circ x $ is a solution of Problem {\em \ref{The Blaschke interpolation problem}}. 		
	\end{enumerate}
\end{thm}
\begin{proof}
By Corollary \ref{CorOfDef10} (3), $|c| \leq |d|$ on $\overline{\mathbb{D}}$. Hence $  \left|\frac{d(\lambda)}{c(\lambda)} \right| \geq 1$ for  $\lambda \in \overline{\mathbb{D}}$. By assumption, $|x_1^{\circ}| < 1$. We claim that $x_1^{\circ} c \neq -d$ on $\overline{\mathbb{D}}$. Suppose that
\begin{align*}
 x_1^{\circ} c = -d \  &\implies  |x_1^{\circ}c| = |d| \\
&\implies   |x_1^{\circ}| |c| = |d| \\
&\implies   |x_1^{\circ}|  = \frac{|d|}{|c|} ,
\end{align*}
which is a contradiction since $\left|\frac{d(\lambda)}{c(\lambda)} \right| \geq 1$ for all $\lambda \in \overline{\mathbb{D}}$, and $|x_1^{\circ}| < 1$ on $\overline{\mathbb{D}}$. Therefore, $x_1^{\circ} c \neq -d$ on $\mathbb{\overline{D}}$.
By Proposition \ref{Prooop}, either $x(\overline{\mathbb{D}}) \subseteq \mathcal{R}_{\overline{\mathcal{E}}}  $ or $x$ is a rational $\overline{\mathcal{{E}}}$-inner function. Because $a,b, c, d$ are polynomials of degree at most $n$, the rational function $x$ has degree at most $n$. 

By the definition of a normalized linear fractional parametrization of the solutions of Problem \ref{The Blaschke interpolation problem}, for some point $\tau \in \mathbb{T} \setminus \{\sigma_1, ... , \sigma_{k}\}$,
$$
\begin{bmatrix}
a(\tau)       & b(\tau) \\
c(\tau)       & d(\tau) \\
\end{bmatrix}
=
\begin{bmatrix}
1 & 0  \\
0 & 1 \\
\end{bmatrix}.
$$
Thus it is easy to see that
\begin{equation} \label{Equ1}
x_1(\tau) = \frac { x_1^{\circ} a(\tau) + b(\tau) } { x_1^{\circ} c(\tau) + d(\tau) }  = x_1^{\circ},
\end{equation}

\begin{equation} \label{Equ2}
x_2(\tau) = \frac{x_3^{\circ} c(\tau) + x_2^{\circ} d(\tau) } { x_1^{\circ} c(\tau) + d(\tau) }  = x_2^{\circ},
\end{equation}

\begin{equation} \label{Equ3}
x_3(\tau) = \frac { x_2^{\circ} b(\tau) + x_3^{\circ} a(\tau)} { x_1^{\circ} c(\tau) + d(\tau) } = x_3^{\circ}.
\end{equation}
By assumption, $|x_3^{\circ}| = 1, |x_1^{\circ}| < 1 $ and $|x_2^{\circ}| < 1$, and hence $x_3(\tau) \neq x_1(\tau) x_2(\tau)$. Therefore, $x(\overline{\mathbb{D}})$ is not in the royal variety $\mathcal{R_{\overline{\mathcal{{E}}}}}$. 

By assumption, $x_2$ is defined by equation (\ref{Equ2}). 
Hence $$ x_2(\sigma_j) = {\displaystyle \frac{x_3^{\circ} c(\sigma_j) + x_2^{\circ} d(\sigma_j)}{x_1^{\circ} c(\sigma_j) + d(\sigma_j)} = \tilde{\eta_j}} \  \text{for} \ j = 1, ... ,n.$$  We want to show that $x$ satisfies the interpolation conditions 
\begin{equation} \label{interCond}
x(\sigma_j) = (\eta_j, \tilde{\eta_j}, \eta_i \tilde{\eta_j}) \quad \mbox{for} \ j = 1, ... ,n, \   
\end{equation}
 which is to say that $\sigma_j, \ j = 1, ... , n$, is a royal node of $x$ with corresponding royal value ($\eta_j, \tilde{\eta_{j}}$). By assumption, there is a finite set $ F \subset \mathbb{T}$ such that, for all $\zeta \in \mathbb{T} \setminus F$, the function
$$
\varphi(\lambda) = 	\psi_\zeta(\lambda) =  \frac { a(\lambda ) \zeta + b(\lambda)} { c(\lambda) \zeta + d(\lambda)}
$$
is a solution of Problem \ref{The Blaschke interpolation problem}, and so
\begin{equation} \label{uppsiz}
\psi_\zeta(\sigma_j) = \eta_j \ \ \text{for} \ j = 1, ..., n
\end{equation}
and
\begin{equation} \label{phasarDD}
A \psi_\zeta(\sigma_j) = \rho_j \ \ \text{for} \ j = 1, ..., k
\end{equation}
for all $\zeta \in \mathbb{T} \setminus F $. Hence, by Proposition \ref{AlgeberaicProp}, 
\begin{equation} \label{bigequation}
\psi_{\zeta(\omega)}(\lambda) = \frac { a(\lambda ) \zeta(\omega) + b(\lambda)} { c(\lambda) \zeta(\omega) + d(\lambda)} = \frac { \omega x_3(\lambda) - x_1(\lambda) }  { \omega x_2(\lambda)   - 1  } = 	\Psi_{\omega} \circ x(\lambda)  
\end{equation}
whenever both sides are defined, that is, for all $\omega \in \mathbb{T} \setminus \Omega_\lambda $ where
$$
\Omega_\lambda = \{ \omega \in \mathbb{T} : \omega x_2({\lambda}) =1 \ \text{or} \ c(\lambda) \zeta(\omega) = -d(\lambda)  \}.
$$	
Note that $\Omega_\lambda$ contains at most two points. Thus, equation (\ref{bigequation}) holds 
as rational functions in $(\omega, \lambda)$, where $\zeta(\omega) = {\displaystyle \frac { \omega x_3^{\circ} - x_1^{\circ} }  { x_2^{\circ}  \omega - 1  }} .$ Hence, for $\omega \in \mathbb{T} \setminus (\zeta^{-1} (F) \cup \Omega_\lambda)$, $	\Psi_{\omega} \circ x $ is a solution of Problem \ref{The Blaschke interpolation problem}, which proves statement (ii).

Equation (\ref{bigequation}) holds for any $\lambda \in \overline{\mathbb{D}}$, that is, provided both denominators are nonzero, and therefore for all but at most two values of $\omega \in \mathbb{T}$.  Combine equations (\ref{uppsiz}) and (\ref{bigequation}) (with $\lambda = \sigma_j$) to infer that, for $j = 1, ... , n$ and for all but finitely many $\omega \in \mathbb{T}$, 
$$
 \frac { \omega x_3(\sigma_j) - x_1(\sigma_j) }  { \omega x_2(\sigma_j)   - 1  } = \psi_{\zeta(\omega)}(\sigma_j) = \eta_j.
$$
Therefore, for almost all $\omega \in \mathbb{T}$ and  $j = 1, ... , n$,
\begin{equation} \label{valueofx's}
\omega x_3(\sigma_j) - x_1(\sigma_j) = \eta_j (\omega x_2(\sigma_j)   - 1).
\end{equation}
Recall that $  x_2(\sigma_j)=  \tilde{\eta_j}  \  \text{for} \ j = 1, ... ,n.$ Hence from equations (\ref{valueofx's}) it follows that $x_1(\sigma_j) = \eta_j$ and $x_3(\sigma_j) = \eta_j \tilde{\eta_j}, \ j=1, ... , n$, and so the interpolation conditions  (\ref{interCond}) hold. 

We have already observed that $x$ is a rational $\overline{\mathcal{{E}}}$-inner function, $\text{deg}(x) \leq n$ and that $x(\overline{\mathbb{D}})$ is not in $\mathcal{R_{\overline{\mathcal{E}}}}$. Thus by Theorem \ref{Degreeofx}, the number of royal nodes of $x$ is equal to the degree of $x$. Consequently $x$ has at most $n$ royal nodes. Because the points $ \sigma_j, j = 1, ... ,n $ are royal nodes, they contain all $n$ royal nodes of $x$, and so deg$(x) =n $. Observe that $k$ of the $\sigma_j$ lie in $\mathbb{T}$; thus $x$ has exactly $k$ royal nodes in $\mathbb{T}.$  Hence $x \in \mathcal{R}^{n,k}$. 

Next we show  that $Ax_1(\sigma_j)     = \rho_j$ for $j = 1, ..., k$. Fix $j \in \{ 1, ..., k\}$. By Proposition \ref{PhasarTheoremForTetra}, for $\omega \in \mathbb{T}, \ \omega \tilde{\eta_j} \neq 1, $ 
\begin{equation} \label{1stEqu}
A(\Psi_{\omega} \circ x) (\sigma_j) = A x_1(\sigma_j).
\end{equation}
In addition  there exists a set $\Omega_j$ containing at most one point in $\T$ such that $c(\sigma_j) \zeta(\omega) + d(\sigma_j) = 0 $ for $\omega \in \Omega_j$. Thus, if $\omega \in \mathbb{T} \setminus ( \{ \overline{{\tilde {\eta}}_j} \} \ \cup \ \Omega_j )$, by equation (\ref{bigequation}),  $\psi_{\zeta(\omega)} = \Psi_{\omega} \circ x $ in a neighbourhood of $\sigma_j$, and so, for such $\omega$, 
\begin{equation} \label{2ndEqu}
A \psi_{\zeta(\omega)}(\sigma_j) = A (\Psi_{\omega} \circ x)(\sigma_j).
\end{equation}
The equations (\ref{1stEqu}), (\ref{2ndEqu}) and (\ref{phasarDD}) all hold for $\omega$ in a cofinite subset of $\mathbb{T}$.   Therefore, for $\omega$ in the intersection of these cofinite subsets,
$$
Ax_1(\sigma_j) = A (\Psi_{\omega} \circ x)(\sigma_j) = A \psi_{\zeta(\omega)} (\sigma_j) = \rho_j.
$$
Thus (i) holds. 
\end{proof}

\begin{cor} \label{LastCor}
	Let $(\sigma, \eta, \rho)$ be Blaschke interpolation data of type $(n,k)$. Let $x$ be a solution of Problem \ref{The royal tetra-interpolation problem} with  data $(\sigma, \eta, \tilde{\eta}, \rho)$ for some $\tilde{\eta_{j}} \in \overline{\mathbb{D}}, j = 1, ..., n, $ and that $x(\overline{\mathbb{D}}) \not\subset \mathcal{R_{\bar{\mathcal{E}}}}$. For all $\omega \in \mathbb{T} \setminus \{ \overline{\tilde{\eta_1}}, ... , \overline{{\tilde {\eta}}_k}\}$, the function $\varphi = \Psi_{\omega} \circ x$ is a solution of Problem \ref{The Blaschke interpolation problem} with Blaschke interpolation data $(\sigma, \eta, \rho)$ . Conversely, for every solution $\varphi$ of the Blaschke interpolation problem with data $(\sigma, \eta, \rho)$, there exists $\omega \in \mathbb{T}$ such that $\varphi = \Psi_{\omega} \circ x$ .
\end{cor}
\begin{proof}
	$(\Longrightarrow)$ Consider Blaschke interpolation data $(\sigma, \eta, \rho)$. If $x= (x_1, x_2, x_3)$ is a solution of Problem \ref{The royal tetra-interpolation problem} with data $(\sigma, \eta, \tilde{\eta}, \rho)$ for some $\tilde{\eta_{j}} \in \overline{\mathbb{D}}, j = 1, ..., n,$ and  $x(\overline{\mathbb{D}}) \not\subseteq \mathcal{R_{\bar{\mathcal{E}}}}$, then, by  Theorem \ref{4.4Tetra} (1), for all $\omega \in \mathbb{T} \setminus \{ \overline{\tilde{\eta_1}}, ... , \overline{{\tilde {\eta}}_k}\}$, there exists a rational function $\varphi = \Psi_{\omega} \circ x$ that solves the Blaschke interpolation problem  with data $(\sigma, \eta, \rho)$. 
	
	$(\Longleftarrow)$ Let $\varphi$ be a solution of the Blaschke interpolation problem (Problem \ref{The Blaschke interpolation problem}) with data $(\sigma, \eta, \rho)$ of type $(n,k)$. Then, by Theorem \ref{4.4Tetra} and Theorem \ref{ReverseThm} (ii), there exists $\omega \in \mathbb{T}$ such that $\varphi = \Psi_{\omega} \circ x$.
\end{proof}

\section{The algorithm}\label{algorithm}

In this section we present  a concrete algorithm for  the solution of the royal $\overline{\mathcal{{E}}}$-interpolation problem. 

Let $(\sigma, \eta, \tilde{\eta}, \rho)$ be royal interpolation data of type $(n,k)$ for the tetrablock, as in Definition \ref{royaltetradata}.  One can consider the associated Blaschke interpolation data $(\sigma, \eta, \rho)$ of type $(n,k)$ . To construct a rational $\overline{\mathcal{{E}}}$-inner function $x: \mathbb{D} \rightarrow \overline{\mathcal{{E}}}$ of degree $n$ having royal nodes $\sigma_j$ for $j = 1, ... , n$, royal values $\eta_{j}, \tilde{\eta_{j}}$, and phasar derivatives $\rho_j$ at $\sigma_j$ for $j = 1, ..., k $, we proceed as follows.

(1) Consider the Pick matrix $M = [m_{i,j}]_{i,j=1}^{n}$ for the data $(\sigma, \eta, \rho)$.  It has entries
\begin{equation}\label{PickMatrix}
m_{i,j} = \begin{cases}
\rho_i & \quad \quad \text{if } i = j \leq k \\
\displaystyle{\frac{1-\overline{\eta_i} \eta_{j}}{1- \overline{\sigma_i} \sigma_j} }& \quad \quad \text{otherwise}. \\
\end{cases}
\end{equation}
Assume that $M$ is positive definite; otherwise the interpolation problem \ref{The Blaschke interpolation problem} is not solvable. Introduce the notation
\begin{equation} \label{xLambdayLambda}
x_\lambda = \begin{bmatrix}
\frac { 1} { 1-\overline{\sigma}_{1}\lambda}  \\ \vdots \\ 
\frac { 1} { 1-\overline{\sigma}_{n}\lambda } \\
\end{bmatrix}, \quad \ \ \ y_\lambda = \begin{bmatrix}
\frac { \overline{\eta}_{1}} { 1-\overline{\sigma}_{1}\lambda }  \\ \vdots \\ 
\frac { \overline{\eta}_{n} } { 1-\overline{\sigma}_{n}\lambda  }
\end{bmatrix}, 
\end{equation}
as in equations (\ref{formula_in3.7}).

(2) Choose a point $\tau \in \mathbb{T} \setminus \{ \sigma_1 , \sigma_2, ..., \sigma_k\}$ such that the set of $\zeta \in \mathbb{T}$ for which 
$$
\langle  M^{-1} x_{\tau}, e_j \rangle = \zeta \langle  M^{-1} y_{\tau}, e_j \rangle \quad \text{for some} \  j \in \{1, ..., n\} 
$$
is finite. Here $e_j$ is the $j$th standard basis vector in $\mathbb{C}^n$. 

(3) Let 
\begin{equation} \label{gEquation}
g(\lambda) = {\displaystyle \prod_{j=1}^{n}} \frac { 1 - \overline{\sigma_j} \lambda} {1 - \overline{\sigma_j} \tau}\ ,
\end{equation}
and let polynomials $a, b, c, d$ be defined by 
\begin{eqnarray} \label{firstpolynomial}
	a(\lambda) &=& g(\lambda) ( 1 - (1 - \overline{\tau} \lambda ) \langle x_{\lambda}, M^{-1} x_{\tau}\rangle), \\ \label{secondpoly}
	b(\lambda) &=& g(\lambda) (1 - \overline{\tau} \lambda ) \langle x_{\lambda}, M^{-1} y_{\tau}\rangle, \\ \label{thirdpoly}
	c(\lambda) &=& -g(\lambda) (1 - \overline{\tau} \lambda ) \langle y_{\lambda}, M^{-1} x_{\tau}\rangle \\ \label{fourpoly}
	d(\lambda) &=& g(\lambda) ( 1 + (1 - \overline{\tau} \lambda ) \langle y_{\lambda}, M^{-1} y_{\tau}\rangle) .
\end{eqnarray} 
Observe that 
\begin{equation} \label{Equationofthm}
\begin{bmatrix}
a(\tau)       & b(\tau) \\
c(\tau)       & d(\tau) \\
\end{bmatrix}
=
\begin{bmatrix}
1 & 0  \\
0 & 1 \\
\end{bmatrix}.
\end{equation}
(See Theorem 3.9 in \cite{ALY3}).

(4) Find $x_1^{\circ}, x_2^{\circ}, x_3^{\circ} \in \mathbb{C}$ such that 
$$
| x_3^{\circ} | = 1, \ \ \ | x_1^{\circ} | < 1, \ \ \  | x_2^{\circ} | < 1, \ \ \ x_1^{\circ} = \overline{x_2^{\circ}} x_3^{\circ},
$$ and 
$$
\frac{x_3^{\circ} c(\sigma_j) + x_2^{\circ} d(\sigma_j)}{x_1^{\circ} c(\sigma_j) + d(\sigma_j)} = \tilde{\eta_j} , \ j = 1, ..., n.  
$$ 
If there is no  $(x_1^{\circ}, x_2^{\circ}, x_3^{\circ})$ satisfying these conditions, then by Theorem \ref{ReverseThm}, the royal $\overline{\mathcal{{E}}}$-interpolation problem is not solvable. 

(5) If there are such $(x_1^{\circ}, x_2^{\circ}, x_3^{\circ}) \in \mathbb{C}$, we define
\begin{eqnarray*}
	x_{1}(\lambda) &=& \frac{x_1^{\circ} a + b}{x_1^{\circ} c + d}(\lambda), \\
	x_{2}(\lambda) &=&  \frac{x_3^{\circ} c + x_2^{\circ} d } { x_1^{\circ} c + d}(\lambda),\\
		x_3(\lambda) &=& \frac { x_2^{\circ} b + x_3^{\circ} a} { x_1^{\circ} c + d }(\lambda), \ \text{for} \  \lambda \in \mathbb{D}.\\ 
\end{eqnarray*}
It is easy to see that, since the equation (\ref{Equationofthm}) is satisfied,
$$
x_1(\tau) = x_1^{\circ}, \  x_2(\tau) = x_2^{\circ} \ \text{and} \ x_3(\tau) = x_3^{\circ}.
$$
Then, by Theorem \ref{ReverseThm},  $x = (x_1, x_2, x_3) $ is a rational $\overline{\mathcal{{E}}}$-inner function of degree at most $n$ such that $x(\sigma_j) = (\eta_j, \tilde{\eta}_j, \eta_j \tilde{\eta}_j)$ for $j = 1, ..., n$, and $Ax_1(\sigma_j) = \rho_j$ for $j = 1, ..., k$. By assumption, $|x_3^{\circ}| = 1, |x_1^{\circ}| < 1 $ and $|x_2^{\circ}| < 1$, and hence $x_3(\tau) \neq x_1(\tau) x_2(\tau)$. Hence, $x(\overline{\mathbb{D}})$ is not in the royal variety of the tetrablock, and the degree of $x$ is exactly $n$. \\

Now let us relate the steps of algorithm to results earlier in the paper. 
\begin{enumerate}
	\item If the royal $\overline{\mathcal{{E}}}$-interpolation problem with data $(\sigma, \eta, \tilde{\eta},  \rho)$ for some $\tilde{\eta_{j}} \in \overline{\mathbb{D}}$ is solvable, then  by Theorem \ref{4.4Tetra}, the Blaschke interpolation problem with data $(\sigma, \eta, \rho)$ is solvable. By \cite[Proposition 3.2] {ALY3}, $M>0$. 
	\item The conditions that 
	$
	| x_3^{\circ} | = 1, \ | x_1^{\circ} | < 1, \  | x_2^{\circ} | < 1,  \text{and} \ x_1^{\circ} = \overline{x_2^{\circ}} x_3^{\circ}
	$ are equivalent to $(x_1^\circ, x_2^\circ, x_3^\circ) \in b\overline{\mathcal{{E}}}$ and $| x_2^{\circ}| < 1 $.  
	\item The equations for $x_1, x_2$ and $x_3$ are equations (\ref{Equationforx_1}), (\ref{Equationforx_2}) and (\ref{Equationforx_3}) respectively.
\end{enumerate}

\section{Examples} \label{2Example}

The following two examples describe all rational tetra-inner functions of degree $1$.

\begin{lem} \label{AutLemma}
	Let $\sigma_{1} \in \mathbb{D} $, and $\eta, \tilde{\eta} \in \overline{\mathbb{D}}$. Let $m \in \mathrm{Aut}(\mathbb{D})$ be such that $m(\sigma_{1}) = 0$. Suppose there exists a rational $\overline{\mathcal{{E}}}$-inner $y : \mathbb{D} \rightarrow \overline{\mathcal{{E}}}$ such that $ y(0) = (\eta, \tilde{\eta}, \eta \tilde{\eta})$. Then $x = y \circ m$ is a rational $\overline{\mathcal{{E}}}$- inner function such that $x(\sigma_{1}) = (\eta, \tilde{\eta}, \eta \tilde{\eta}) $. 
\end{lem}
	\begin{proof}
		By assumption,  the function $y : \mathbb{D} \rightarrow \overline{\mathcal{{E}}} $ is such that $y(0) = (\eta, \tilde{\eta}, \eta \tilde{\eta} )$.   The Blaschke factor $m : \mathbb{D} \rightarrow \mathbb{D}$ such that $m(z) = \displaystyle{\frac{z- \sigma_{1}}{1 - \overline{\sigma_{1}} z}}$ moves $\sigma_{1}$  to $0$.
		\newline Note that $(y \circ m ) (\sigma_{1}) = y(m(\sigma_{1})) = y(0) = (\eta, \tilde{\eta}, \eta \tilde{\eta}). $
		It is easy to see that the composition $x = y \circ m$ is a rational $\overline{\mathcal{{E}}}$- inner function,   $x : \mathbb{D} \rightarrow \overline{\mathcal{{E}}}$ such that $x(\sigma_{1}) = (\eta, \tilde{\eta}, \eta \tilde{\eta})$. 
	\end{proof}

\begin{exmp} \label{FirstExample}
	{\em Let $(\sigma, \eta, \tilde{\eta}, \rho)$ be royal interpolation data of type $(1,0)$ for the tetrablock, as in Definition \ref{royaltetradata}.
Consider  Problem \ref{The royal tetra-interpolation problem} with the data $(\sigma, \eta, \tilde{\eta}, \rho)$. We are given  a single royal node $\sigma_{1} \in \mathbb{D}$ and a  royal value  $ (\eta, \tilde{\eta}, \eta \tilde{\eta})$, where $\eta, \tilde{\eta} \in \mathbb{D}$, and we need to find a $\overline{\mathcal{{E}}}$-inner function $x$ of degree $1$ such that $x(\sigma_{1}) = (\eta, \tilde{\eta}, \eta \tilde{\eta})$. By composition with an automorphism of $\mathbb{D}$, we may reduce our problem  to the case that $\sigma_{1} = 0$.
	
{\bf Step 1}. Pick an arbitrary $\tau \in \mathbb{T}$. The  solution set of the associated Blaschke interpolation problem has the normalized parametrization, given, for some $\zeta \in \mathbb{T}$, by
$$
\varphi(\lambda) = \frac { a(\lambda)\zeta + b(\lambda)} { c(\lambda)\zeta + d(\lambda)}, \quad \text{for all} \ \ \lambda \in \mathbb{D},
$$
where  $a, b, c, d$ are given by equations  \eqref{firstpolynomial} - \eqref{fourpoly},
 and $x_\lambda, y_\lambda, g$ and $M$ are given by equations (\ref{xLambdayLambda}), (\ref{gEquation}) and (\ref{PickMatrix}) respectively.  
Note that since  $\sigma_1 = 0$, $g(\lambda) = {\displaystyle \frac { 1 - \overline{\sigma_1} \lambda} {1 - \overline{\sigma_1} \tau}}\ = 1$,  $M= \displaystyle{\frac{1-\overline{\eta} \eta}{1- \overline{\sigma_1} \sigma_1} }= 1-|\eta|^{2}$, $x_\lambda = \displaystyle{\frac{1}{1 - 0 \lambda} = 1 }$ and $y_\lambda = \displaystyle{\frac{\overline{\eta}}{1 - 0 \lambda} = \overline{\eta}} $. 

Thus, the polynomials $a, b, c$ and $d$ are defined by 
\begin{eqnarray}\label{Valuofa}
a(\lambda)\nonumber &=& g(\lambda) ( 1 - (1 - \overline{\tau} \lambda ) \langle x_{\lambda}, M^{-1} x_{\tau}\rangle) \\ 
&=& 1 - \frac{1-\overline{\tau} \lambda}{1 - |\eta|^{2}} 
= \frac{\overline{\tau} \lambda - |\eta|^{2}}{1 - |\eta|^{2}}, 
\end{eqnarray}
\begin{eqnarray}\label{Valuofb}
b(\lambda)\nonumber &=& g(\lambda) (1 - \overline{\tau} \lambda ) \langle x_{\lambda}, M^{-1} y_{\tau}\rangle \nonumber\\
&=& 1 (1 - \overline{\tau} \lambda) \frac{\eta}{1- |\eta|^{2}} 
= \frac{\eta (1 - \overline{\tau} \lambda)}{1 - |\eta|^{2}}, 
\end{eqnarray}
\begin{eqnarray}\label{Valuofc}
c(\lambda)\nonumber &=& -g(\lambda) (1 - \overline{\tau} \lambda ) \langle y_{\lambda}, M^{-1} x_{\tau}\rangle \nonumber\\
&=& -1 (1 - \overline{\tau} \lambda) \frac{\overline{\eta}}{1- |\eta|^{2}} 
= \frac{-\overline{\eta} (1 - \overline{\tau} \lambda)}{1 - |\eta|^{2}}, 
\end{eqnarray}
\begin{eqnarray}\label{Valuofd}
d(\lambda)\nonumber &=& g(\lambda) ( 1 + (1 - \overline{\tau} \lambda ) \langle y_{\lambda}, M^{-1} y_{\tau}\rangle) \\
&=& 1 ( 1 + (1 - \overline{\tau} \lambda) \frac{\eta \overline{\eta}}{1- |\eta|^{2}}) =  \frac{1 - |\eta|^{2} \overline{\tau} \lambda}{1 - |\eta|^{2}}.
\end{eqnarray}
{\bf Step 2}.  The next step is to determine whether there exist $x_1^{\circ}, x_2^{\circ}, x_3^{\circ} \in \mathbb{C}$ such that 
\begin{equation}\label{x_zero_1}
| x_3^{\circ} | = 1, \ \ \ | x_1^{\circ} | < 1, \ \ \  | x_2^{\circ} | < 1, \  \ x_1^{\circ} = \overline{x_2^{\circ}} x_3^{\circ},
\end{equation} and 
\begin{equation}\label{x_zero_2}
\frac{x_3^{\circ} c(0) + x_2^{\circ} d(0)}{x_1^{\circ} c(0) + d(0)} = \tilde{\eta} .  
\end{equation}
Here, $$a(0) = \displaystyle{\frac{- |\eta|^2}{1- |\eta|^2}}, \quad  b(0) = \displaystyle{\frac{ \eta}{1- |\eta|^2}}, \quad  c(0) = \displaystyle{\frac{- \overline{\eta}}{1- |\eta|^2}}, \quad  d(0) = \displaystyle{\frac{1}{1- |\eta|^2}}.$$
Let $x_3^{\circ} = \omega$ for $\omega \in \mathbb{T}$. Thus
\begin{eqnarray} \nonumber
\frac{x_3^{\circ} c(0) + x_2^{\circ} d(0)}{x_1^{\circ} c(0) + d(0)} = \tilde{\eta} \
&\Leftrightarrow& x_2^{\circ} = - x_1^{\circ} \overline{\eta} \tilde{\eta} + \tilde{\eta} + \omega \overline{\eta}. \\\nonumber
\end{eqnarray}
Since $x_1^{\circ} = \overline{x_2^{\circ}} x_3^{\circ}$ and $x_3^{\circ} = \omega$,  we have the system
\begin{equation}
\left\{
\begin{array}{lcl}
x_3^{\circ}  & = & \omega \\
x_1^{\circ} & = & \overline{x_2^{\circ}} \omega \\
x_2^{\circ}    & = & - x_1^{\circ} \overline{\eta} \tilde{\eta} + \tilde{\eta} + \omega \overline{\eta}.
\end{array}
\right.
\end{equation}
For given $\eta, \tilde{\eta} \in \mathbb{D}$, we want to find a solution $x_1^{\circ}, x_2^{\circ}, x_3^{\circ}$ of the above system  such that $x_3^{\circ} \in \mathbb{T}$, $|x_1^{\circ}| <1$, and $|x_2^{\circ}| < 1$.
This is equivalent to finding $\omega \in \mathbb{T}$ and $|x_1^{\circ}| < 1$, such that the equation
\begin{equation} \label{SatisfyingEquu}
\omega \overline{\tilde{\eta}} + \eta= x_1^{\circ} + \overline{ x_1^{\circ}} \omega   \overline{\tilde{\eta}} \eta 
\end{equation}
holds.
\begin{lemma} \label{beta}
Let $\eta, \tilde{\eta} \in \mathbb{D}$. There are $\omega \in \mathbb{T}$ and $|x_1^{\circ}| < 1$ such that equation \eqref{SatisfyingEquu} holds.
\end{lemma}
\begin{proof}
Choose $\omega \in \mathbb{T}$. 
Let $\xi= \omega \overline{\tilde{\eta}}$, and so $\xi \in \mathbb{D}$.
Let $s= \xi +\eta$ and $p=\xi \eta$; then $(s, p)$ is a point of the symmetrized bidisc $G= \{(z_1 +z_2, z_1 z_2): |z_1| < 1,  |z_2| < 1\}$. One can see that
$$\beta = \frac{s- \bar{s}p}{1 - |p|^2}$$
satisfies the equation $ s =\beta+ \bar{\beta} p$ and the inequality $|\beta| <1$. Hence,
for $\omega \in \mathbb{T}$ and for $x_1^{\circ}= \beta$,
equation \eqref{SatisfyingEquu} holds.
\end{proof}

{\bf Step 3}. Therefore, for the given data, $0 \rightarrow (\eta, \tilde{\eta}, \eta \tilde{\eta}) $, there is a $1$-parameter family of  $(x_1^{\circ}, x_2^{\circ},x_3^{\circ})$ such that equations \eqref{x_zero_1} and \eqref{x_zero_2} satisfied, given by
$$x_1^{\circ}=\frac{(\omega \overline{\tilde{\eta}}+ \eta) - ( \bar{\omega} 
\tilde{\eta}+ \bar{\eta})\omega \overline{\tilde{\eta}}\eta}{1 - |\tilde{\eta} \eta|^2},
\ x_2^{\circ}=\omega  \overline{x_1^{\circ}}, \ x_3^{\circ}=\omega,
$$ 
for any $\omega \in \mathbb{T}$.
Substitute these values into equations  (\ref{Equationforx_1}), (\ref{Equationforx_2}) and (\ref{Equationforx_3}) to obtain the degree $1$ rational 
$\overline{\mathcal{{E}}}$-inner function $x= (x_1, x_2, x_3)$, satisfying $x(0)=(\eta, \tilde{\eta}, \eta \tilde{\eta})$, where, for $\lambda \in \mathbb{D}$,
\begin{eqnarray}\label{sol-example1}
x_1(\lambda) &=& \frac{x_1^{\circ} (\overline{\tau} \lambda - |\eta|^{2})+ \eta (1 - \overline{\tau} \lambda) }{-x_1^{\circ}\overline{\eta} (1 - \overline{\tau} \lambda) +  1 - |\eta|^{2} \overline{\tau} \lambda }, \\
x_2(\lambda) &=&
\frac{-x_3^{\circ}\overline{\eta} (1-\overline{\tau} \lambda) + x_2^{\circ} ( 1- |\eta|^{2} \overline{\tau} \lambda)}{-x_1^{\circ} \overline{\eta} (1-\overline{\tau} \lambda) + 1- |\eta|^{2} \overline{\tau} \lambda },\\
x_3(\lambda) &=&
\frac{x_2^{\circ}\eta (1-\overline{\tau} \lambda) + x_3^{\circ} (\overline{\tau} \lambda - |\eta|^{2})}{-x_1^{\circ} \overline{\eta} (1 - \overline{\tau} \lambda) + 1 - |\eta|^{2}\overline{\tau} \lambda  }.
\end{eqnarray}
}
\end{exmp}

\begin{exmp} \label{SecondExample}
	{\em Consider the case $n=1, k=1$. Suppose $\sigma = 1$. The points $\eta, \tilde{\eta}  \in \mathbb{T}$ and  $\rho > 0 $ are prescribed, and we seek a $\overline{\mathcal{{E}}}$-inner function $x =(x_1, x_2, x_3)$ of degree $1$ such that $x(1) = (\eta, \tilde{\eta}, \eta \tilde{\eta})$ and $Ax_1(1) = \rho$. 
		
{\bf Step 1.} Choose $\tau \in \mathbb{T} \ \setminus \{1\} $. The normalized parametrization of the solution set of the associated Blaschke interpolation problem is given by
		\begin{equation}\label{phiEquationn}
		\varphi(\lambda) = \frac { a(\lambda)\zeta + b(\lambda)} { c(\lambda)\zeta + d(\lambda)} \quad \mbox{for} \ \lambda \in \mathbb{D}, \mbox{and some} \ \zeta \in \mathbb{T},
		\end{equation}
where  $a, b, c, d$ are given by equations  \eqref{firstpolynomial} - \eqref{fourpoly},
 and $x_\lambda, y_\lambda, g$ and $M$ are given by equations (\ref{xLambdayLambda}), (\ref{gEquation}) and (\ref{PickMatrix}) respectively.  

Note that, since  $\sigma = 1$ and $ k = 1$, $g(\lambda) = {\displaystyle  \frac { 1 - \overline{\sigma} \lambda} {1 - \overline{\sigma} \tau}}\ = {\displaystyle\frac{1-\lambda}{1-\tau}}$, $M= \rho $,  $x_\lambda = \displaystyle{\frac{1}{1-\lambda}}$ and $y_\lambda = \displaystyle{\frac{\overline{\eta}}{1-\lambda}}$.
 Therefore, polynomials $a, b, c$ and $d$ are defined by
		\begin{eqnarray} \label{secondA}
			a(\lambda) &=& g(\lambda) \Big( 1 - (1 - \overline{\tau} \lambda ) \langle x_\lambda, M^{-1} x_\tau \rangle \Big) \nonumber\\
			&=&  \frac{1-\lambda}{1-\tau} - \frac{ (1 - \overline{\tau} \lambda)}{\rho |1-\tau|^{2} }    
		\end{eqnarray}
	\begin{eqnarray}\label{secondB}
		b(\lambda) &=& g(\lambda) (1 - \overline{\tau} \lambda ) \langle x_{\lambda}, M^{-1} y_{\tau}\rangle \nonumber \\
		&=&  \frac{\eta(1-\overline{\tau} \lambda)}{\rho |1-\tau|^{2}}     
	\end{eqnarray}
\begin{eqnarray}\label{secondC}
	c(\lambda) &=& -g(\lambda) (1 - \overline{\tau} \lambda ) \langle y_{\lambda}, M^{-1} x_{\tau}\rangle \nonumber \\
	&=&  - \frac{\overline{\eta} (1- \overline{\tau} \lambda)}{\rho |1-\tau|^{2}}  
\end{eqnarray}
\begin{eqnarray}\label{secondD}
	d(\lambda) &=& g(\lambda) \Big( 1 + (1 - \overline{\tau} \lambda ) \langle y_{\lambda}, M^{-1} y_{\tau}\rangle \Big) \nonumber \\
	&=&  \frac{1-\lambda}{1-\tau} + \frac{1 - \overline{\tau} \lambda}{\rho |1-\tau|^{2}}.   
\end{eqnarray}
{\bf Step 2.} Next, determine if there exist $x_1^{\circ}, x_2^{\circ}, x_3^{\circ} \in \mathbb{C}$ such that 
\begin{equation}\label{x_zero_3}
| x_3^{\circ} | = 1, \ \ \ | x_1^{\circ} | < 1, \ \ \  | x_2^{\circ} | < 1, \ \text{and} \ x_1^{\circ} = \overline{x_2^{\circ}} x_3^{\circ},
\end{equation}
and 
\begin{equation}\label{x_zero_4}
\frac{x_3^{\circ} c(1) + x_2^{\circ} d(1)}{x_1^{\circ} c(1) + d(1)} = \tilde{\eta} .  
\end{equation}
Here, $$a(1) = \displaystyle{\frac{1- \overline{\tau}}{\rho |1-\overline{\tau}|^{2}}}, \quad  b(1) = \displaystyle{\frac{ \eta - \eta \overline{\tau}}{\rho |1-\overline{\tau}|^{2}}}, \quad \ c(1) = \displaystyle{\frac{- \overline{\eta} + \overline{\eta} \overline{\tau} }{\rho |1-\overline{\tau}|^{2}}}, \quad  d(1) = \displaystyle{\frac{1- \overline{\tau}}{\rho |1-\overline{\tau}|^{2}}}.$$
Let $x_3^{\circ} = \omega$ for $\omega \in \mathbb{T}$. Thus
\begin{eqnarray} \nonumber
\frac{x_3^{\circ} c(1) + x_2^{\circ} d(1)}{x_1^{\circ} c(1) + d(1)} = \tilde{\eta} \ 	&\Leftrightarrow& 
\ x_2^{\circ} =  -\overline{\eta} x_1^{\circ} \tilde{\eta}  + \tilde{\eta}+ \overline{\eta} \omega.   \\\nonumber
\end{eqnarray}
Since $x_1^{\circ} = \overline{x_2^{\circ}} x_3^{\circ}$ and $x_3^{\circ} = \omega$,  we have 
\begin{equation} \label{SystemExample4}
\left\{
\begin{array}{lcl}
x_3^{\circ}  & = & \omega \\
x_1^{\circ} & = & \overline{x_2^{\circ}} \omega \\
x_2^{\circ}    & = &  -\overline{\eta} x_1^{\circ} \tilde{\eta}  + \tilde{\eta}+ \overline{\eta} \omega.
\end{array}
\right.
\end{equation}
For given $\eta, \tilde{\eta} \in \mathbb{T}$, we want to find a solution $x_1^{\circ}, x_2^{\circ}, x_3^{\circ}$ of the above system  such that $|x_1^{\circ}| <1, |x_2^{\circ}| < 1$. Thus we want  to find $\omega \in \mathbb{T}$ and $x_1^{\circ} \in \mathbb{C} $ such that $|x_1^{\circ}|<1$ and 
\begin{equation} \label{SatisEq22}
x_1^{\circ} + \overline{x_1^{\circ}} \omega  \overline{\tilde{\eta}} \eta = \omega \overline{\tilde{\eta}} + \eta.
\end{equation}
\begin{lemma} \label{beta-2}
Let $\eta, \tilde{\eta} \in \mathbb{T}$. There are $\omega \in \mathbb{T}$ and $x_1^{\circ} \in \D$ such that equation \eqref{SatisEq22} holds.
\end{lemma}
\begin{proof}
Choose $\omega \in \mathbb{T}$ such that $\omega \neq \eta \tilde{\eta}$. 
Let $\xi= \omega \overline{\tilde{\eta}}$, and so $\xi \in \mathbb{T}$ and $-\xi \bar{\eta} \neq -1$.
 We seek $\beta \in \D$ such that 
\begin{equation}\label{eq-for-beta}
\beta + \bar{\beta} \xi \eta =\xi +\eta;
\end{equation}
that is, we seek $\beta \in \D$ such that the following equation is satisfied
\begin{equation}\label{eq-for-beta-2}
- \xi \bar{\eta} = \frac{ 1- \beta \bar{\eta}} {1 - \bar{\beta} \eta}.
\end{equation}
Let $- \xi \bar{\eta} = e^{i \phi}$ where $ 0 \le \phi < 2\pi$ and $\phi \neq \pi$. Hence we can choose $\beta$ such that 
$ 1- \beta \bar{\eta} = r  e^{i \phi/2}$, that is, 
$\beta = \eta(1-r  e^{i \phi/2})$.
Choose $r$ so small that $|\beta| <1$. 
 In the case $ 0 \le \phi/2 < \pi/2$,
one can find $r>0$ so small that $|\beta|<1$.
In the case $ \pi/2 < \phi/2 < \pi$,
one can find $r<0$ so small that $|\beta| <1$.

Therefore, for $\omega \in \mathbb{T}$  such that $\omega \neq \eta \tilde{\eta}$ and for $x_1^{\circ}= \beta = \eta (1-r e^{i \phi/2})$,
equation \eqref{eq-for-beta} holds, and so the equivalent equation \eqref{SatisEq22} holds too.
\end{proof}

{\bf Step 3}. Therefore, for the given data, $1 \rightarrow (\eta, \tilde{\eta}, \eta \tilde{\eta}) $ and $ \rho >0$, there is a family of points  $(x_1^{\circ}, x_2^{\circ},x_3^{\circ})$ such that equations \eqref{x_zero_3} and \eqref{x_zero_4} are satisfied. 
Substitution of these values into equations  (\ref{Equationforx_1}), (\ref{Equationforx_2}) and (\ref{Equationforx_3}) yields the degree $1$ rational 
$\overline{\mathcal{{E}}}$-inner function $x= (x_1, x_2, x_3)$ where, for $\lambda \in \mathbb{D}$,
\begin{eqnarray} \label{Firestsubtttt}
x_{1}(\lambda) &=& \frac{x_1^{\circ}[(1-\lambda)\rho (1-\overline{\tau})-(1-\overline{\tau} \lambda)] + \eta(1-\overline{\tau} \lambda)}{x_1^{\circ} [- \overline{\eta} (1-\overline{\tau} \lambda) ] + \rho (1-\overline{\tau}) (1-\lambda) +(1-\overline{\tau} \lambda) }, 
\end{eqnarray}
\begin{eqnarray}\label{Secondsubttttt}
x_{2}(\lambda) &=&   \frac{x_3^{\circ}[- \overline{\eta} (1-\overline{\tau} \lambda)] + x_2^{\circ} [\rho (1-\overline{\tau}) (1-\lambda) +(1-\overline{\tau} \lambda)]}{x_1^{\circ} [- \overline{\eta} (1-\overline{\tau} \lambda) ] + \rho (1-\overline{\tau}) (1-\lambda) +(1-\overline{\tau} \lambda) },
\end{eqnarray}
\begin{eqnarray}\label{Thirdsubttttt}
x_3(\lambda) &=&
\frac{x_2^{\circ} \eta (1-\overline{\tau} \lambda) + x_3^{\circ} [\rho (1-\overline{\tau}) (1-\lambda) -(1-\overline{\tau} \lambda)]}{x_1^{\circ} [- \overline{\eta} (1-\overline{\tau} \lambda) ] + \rho (1-\overline{\tau}) (1-\lambda) +(1-\overline{\tau} \lambda) }.
\end{eqnarray}
One can  check that $x= (x_1,x_2,x_3)$ defined by equations (\ref{Firestsubtttt}), (\ref{Secondsubttttt}) and (\ref{Thirdsubttttt})
is a $\overline{\mathcal{{E}}}$-inner function of degree $1$ satisfying $x(1)=(\eta, \tilde{\eta}, \eta \tilde{\eta})$ and $Ax_1(1) = \rho$.

}\end{exmp}

\bibliography{references}

HADI O. ALSHAMMARI, Jouf University, King Khaled Road, Skaka,
Kingdom of Saudi Arabia; e-mail: hahammari@ju.edu.sa\\

ZINAIDA A. LYKOVA,
School of Mathematics, Statistics and Physics, Newcastle University, Newcastle upon Tyne
 NE\textup{1} \textup{7}RU, U.K.~~\\

\end{document}